\documentclass[12pt]{amsart}
\usepackage{amssymb,latexsym}
\usepackage{multicol,ragged2e}
\usepackage{graphicx,color}

\usepackage{epic,eepic}

\newcommand{\light}[1]{{#1}}

\definecolor{light}{rgb}{1,0,0}        
\renewcommand{\light}[1]{\color{light}{#1}\color{black}}

\numberwithin{equation}{section}

\newtheorem{theorem}{Theorem}[section]
\newtheorem{proposition}[theorem]{Proposition}

\newtheorem{corollary}[theorem]{Corollary}
\newtheorem{lemma}[theorem]{Lemma}

\theoremstyle{definition}

\newtheorem{remark}[theorem]{Remark}
\newtheorem{example}[theorem]{Example}
\newtheorem{definition}[theorem]{Definition}

\def\ZZ{\mathbb{Z}}
\def\CC{\mathbb{C}}

\def\QQ{\mathbb{Q}}
\def\PP{\mathbb{P}}

\def\lg{\mathfrak{g}}

\def\hh{\mathfrak{h}}

\def\dd{\mathbf{d}}

\def\gg{\mathbf{g}}
\def\xx{\mathbf{x}}
\def\yy{\mathbf{y}}
\def\ii{\mathbf{i}}

\def\Acal{\mathcal{A}}
\def\Fcal{\mathcal{F}}
\def\Xcal{\mathcal{X}}

\def\Trop{\operatorname{Trop}}

\def\Poly{\mathbf{P}}

\def\l{\ell}
\def\wnot{w_\circ}

\newcommand{\beal}{\begin{eqnarray}\begin{array}{l} }
\newcommand{\bear}{\begin{eqnarray}\begin{array}{r} }
\newcommand{\beac}{\begin{eqnarray}\begin{array}{c} }
\newcommand{\bealn}{\[\begin{array}{l} }
\newcommand{\bearn}{\[\begin{array}{r} }
\newcommand{\beacn}{\[\begin{array}{c} }
\newcommand{\eea}{\end{array}\end{eqnarray}}
\newcommand{\eean}{\end{array}\]}

\newcommand{\beq}{\begin{equation} }
\newcommand{\eeq}{\end{equation} }

\newcommand{\mat}[4]{\left(\!\!\begin{array}{cc}
#1 & #2 \\
#3 & #4 \\
\end{array}\!\!\right)}

\renewcommand{\eqref}[1]{{\rm (\ref{#1})}}

\begin{document}

\title[Cluster algebras of finite type via Coxeter elements]
{Cluster algebras of finite type via Coxeter elements and principal minors}

\author{Shih-Wei Yang}
\address{\noindent Department of Mathematics, Northeastern University,
 Boston, MA 02115}
\email{yang.s@neu.edu}

\author{Andrei Zelevinsky}
\address{\noindent Department of Mathematics, Northeastern University,
 Boston, MA 02115}
\email{andrei@neu.edu}

\begin{abstract}
We give a uniform geometric realization for
the cluster algebra of an arbitrary finite type with principal
coefficients at an arbitrary acyclic seed.
This algebra is realized as the coordinate ring of a certain reduced double Bruhat cell
in the simply connected semisimple algebraic group of the same Cartan-Killing type.
In this realization, the cluster variables appear as
certain (generalized) principal minors.
\end{abstract}


\date{April 20, 2008; Revised May 19, 2008}

\dedicatory{To Bertram Kostant on the occasion of his 80th birthday}

 \thanks{Both authors supported by A.~Zelevinsky's NSF (DMS) grant \#
 0500534; A.~Zelevinsky also supported by a Humboldt Research Award.}

\maketitle


\tableofcontents


\section{Introduction and main results}
\label{sec:introduction}

The theory of cluster algebras has a lot in common with the theory
of Kac-Moody algebras. In both instances, the structure of an
algebra in question is encoded by a square integer matrix: a
generalized Cartan matrix~$A$ in the Kac-Moody case, and an
exchange matrix~$B$ in the case of cluster algebras.
(Note an important distinction between the two cases: the
sign pattern of matrix entries is symmetric for $A$ but
skew-symmetric for $B$.) In both instances, there is a natural
notion of finite type, and the algebras of finite type share the
same Cartan-Killing classification. This occasionally leads to an
intriguing collision of two types of symmetry: the coordinate ring of a variety
associated with a semisimple algebraic group~$G$ may carry a natural cluster algebra structure
whose Cartan-Killing type is completely different from that
of~$G$. For example, the base affine space $G/N$ for $G = SL_5$
inherits the symmetry of type~$A_4$ from~$G$, but the coordinate
ring $\CC[G/N]$ has a natural cluster algebra structure of
type~$D_6$, see \cite[Proposition~2.26]{ca3}.

One way to relate the two classifications to each other is
by exhibiting a uniform construction of a variety associated with a
semisimple group~$G$ whose coordinate ring naturally carries a cluster algebra
structure of the same Cartan-Killing type.
Such a construction was given in \cite[Example~2.24]{ca3} with
the variety in question being the double
Bruhat cell $G^{c,c^{-1}}$, where $c$ is a Coxeter element  in the
Weyl group of $G$ (here and in the sequel, we assume $G$ to be simply connected).
This result serves as a point of departure for the current paper.
Here we improve on it in the following two aspects.

First, there are many non-isomorphic cluster algebras of the same
type differing from each other by the choice of a coefficient system.
From this perspective, there is nothing especially distinguished
about the coefficient system for $\CC[G^{c,c^{-1}}]$.
On the other hand, as shown in \cite{ca4}, with any exchange matrix
one can associate an important system of \emph{principal}
coefficients which in a certain sense controls all the other
possible choices of coefficients.
In this paper we show that one can realize the cluster algebra of finite type
with principal coefficients at an arbitrary acyclic initial
cluster by replacing the double cell $G^{c,c^{-1}}$ with its
\emph{reduced} version $L^{c,c^{-1}}$ introduced in \cite{bz01}.

Second and perhaps more importantly, the treatment in \cite{ca3}
was limited to the description of the initial cluster in $\CC[G^{c,c^{-1}}]$,
with no information provided about the rest of the cluster variables.
Here we compute explicitly all cluster variables in $\CC[L^{c,c^{-1}}]$.
They turn out to be an interesting special family of
principal generalized minors.

Before stating the general results, we present a motivating example.
Let $G = SL_{n+1}(\CC)$ be of type $A_n$, and let $c = s_1 \cdots s_n$,
in the standard numbering of simple roots.
In this case $L^{c,c^{-1}}$ is
the subvariety of $SL_{n+1}(\CC)$ consisting of tridiagonal
matrices of the form
\begin{equation}
\label{eq:tridiagonal}
M = \begin{pmatrix}v_1 & y_1 & 0 & \cdots & 0 \\
                  1  &  v_2 & y_2 & \ddots & \vdots \\
                  0  & 1 & \ddots & \ddots  & 0 \\
                  \vdots & \ddots & \ddots & \ddots  & y_n\\
                  0 &  \cdots &  0 & 1 & v_{n+1}
                  \end{pmatrix}
\end{equation}
with all $y_1, \dots, y_n$ non-zero.
Let $x_{[i,j]} \in \CC[L^{c,c^{-1}}]$ denote the regular function
on $L^{c,c^{-1}}$ given by the minor with rows and columns $i,
i+1, \dots, j$, with the convention that $x_{[i,j]} = 1$ unless
$1 \leq i \leq j \leq n+1$.
For instance, we have $x_{[i,i]} = v_i$, and $x_{[1,n+1]} = \det(M) = 1$.

In the following result we use the terminology on cluster algebras
introduced in \cite{ca4}; the relevant notions will be recalled
in Section~\ref{sec:Cluster-algebras}.

\begin{theorem}\
\label{th:ca-principal-minors-A-special}
\begin{enumerate}
\item The algebra $\Acal = \CC[L^{c,c^{-1}}]$ is the cluster algebra of type $A_n$ with
principal coefficients at the initial seed $(\xx, \yy, B)$
given by
$$\xx = (x_{[1,1]}, x_{[1,2]}, \dots, x_{[1,n]}),$$

$$\yy = (y_1, \dots, y_n),$$

$$B = \begin{pmatrix}0 & 1 & 0 & \cdots & 0\,\, \\
                  -1  &  0 & 1 & \ddots & \vdots\,\, \\
                  0  & -1 & \ddots & \ddots  & 0\,\, \\
                  \vdots & \ddots & \ddots & \ddots  & 1\,\,\\
                  0 &  \cdots &  0 & -1 & 0\,\,
                  \end{pmatrix}.$$
\item The set of cluster variables in $\Acal$ is $\{x_{[i,j]}: 1 \leq i
\leq j \leq n+1, \ (i,j) \neq (1,n+1)\}$.
\item The exchange relations in $\Acal$ are:
\begin{equation}
\label{eq:exchange-A-principal}
x_{[i,k]} x_{[j,\ell]} = y_{j-1} y_j \cdots y_k x_{[i,j-2]} x_{[k+2,\ell]} + x_{[i,\ell]} x_{[j,k]}
\end{equation}
for $1 \leq i \leq j-1 \leq k \leq \ell - 1 \leq n$.
\end{enumerate}
\end{theorem}

Note that among the relations \eqref{eq:exchange-A-principal},
there are the exchange relations from the initial cluster, which
can be rewritten as follows:
\begin{equation}
\label{eq:exchange-A-initial}
x_{[1,k+1]}  = v_{k+1} x_{[1,k]}  - y_k x_{[1,k-1]}  \quad (k = 1,
\dots, n).
\end{equation}
These relations play a fundamental part in the classical theory
of orthogonal polynomials in one variable.
Thus, the cluster algebra in Theorem~\ref{th:ca-principal-minors-A-special}
can be viewed as some kind of ``enveloping algebra" for this theory.
Note also that the variety $L^{c,c^{-1}}$ is very close to the
variety of tridiagonal matrices used by B.~Kostant \cite{kostant-2}
in his study of a generalized Toda lattice.

Our main result is a generalization of Theorem~\ref{th:ca-principal-minors-A-special}
to the case of an arbitrary simply connected semisimple complex
group~$G$ and an arbitrary Coxeter element~$c$ in its Weyl group~$W$.
We use the terminology and notation in \cite[Sections~4.2, 4.3]{bz01}
(more details are given in Sections~\ref{sec:cce} and \ref{sec:primitive-exchange-relations}).
Recall that~$W$ is a finite Coxeter group generated by the simple
reflections $s_i \, (i \in I)$, and that $c = s_{i_1} \cdots
s_{i_n}$ for some permutation $(i_1, \dots, i_n)$ of the index set~$I$.
For $i, j \in I$, we will write $i \prec_c j$ if $i$ and $j$ are joined by an edge
in the Coxeter graph (that is, if the Cartan matrix entry $a_{i,j}$ is non-zero), and $s_i$
precedes $s_j$ in the factorization of~$c$.
We associate to~$c$ a skew-symmetrizable matrix $B(c) = (b_{i,j})_{i,j \in I}$
by setting
\begin{equation}
\label{eq:Bc}
b_{i,j} = \begin{cases}
-a_{i,j} & \text{if $i \prec_c j$;} \\[.05in]
a_{i,j} & \text{if $j \prec_c i$;} \\[.05in]
0 & \text{otherwise.}
\end{cases}
\end{equation}

Recall from \cite{fz-double,bz01} the finite family of (generalized)
minors $\Delta_{\gamma, \delta} \in \CC[G]$ labeled by pairs of
weights $\gamma, \delta$ belonging to the $W$-orbit of the same
fundamental weight $\omega_i \, (i \in I)$.
We call a minor $\Delta_{\gamma, \delta}$ \emph{principal} if
$\gamma = \delta$; note that this terminology differs from that in
\cite{fz-double}, where ``principal" referred only to
the minors $\Delta_{\omega_i, \omega_i}$.
We will use the following notation:
\begin{equation}
\label{eq:x-gamma-c}
\text{$x_{\gamma;c} \in \CC[L^{c,c^{-1}}]$ is the restriction
of $\Delta_{\gamma, \gamma}$ to the reduced double cell $L^{c,c^{-1}}$.}
\end{equation}
As shown in \cite{fz-double,bz01}, each of the minors
$\Delta_{\omega_i, c \omega_i}$ vanishes nowhere on
$L^{c,c^{-1}}$, so can be viewed as an invertible element
of $\CC[L^{c,c^{-1}}]$.
For $j \in I$, we use the following notation:
\begin{equation}
\label{eq:y-j-c}
\text{$y_{j;c} \in \CC[L^{c,c^{-1}}]$ is
the restriction to  $L^{c,c^{-1}}$ of the product
$\Delta_{\omega_j, c \omega_j} \prod_{i \prec_c j}
\Delta_{\omega_i, c \omega_i}^{a_{i,j}}$.}
\end{equation}

Now we are ready to generalize Part (1) of
Theorem~\ref{th:ca-principal-minors-A-special}.

\begin{theorem}
\label{th:ca-general-1}
The coordinate ring $\Acal(c) = \CC[L^{c,c^{-1}}]$ is the cluster algebra with
principal coefficients at the initial seed $(\xx, \yy, B)$
given by
$$\xx = (x_{\omega_i;c}: i \in I), \quad \yy = (y_{i;c}: i \in I),
\quad B = B(c).$$
Furthermore, $\Acal(c)$ is of finite type, and its Cartan-Killing
type is the same as that of~$G$.
\end{theorem}

It is easy to see that every acyclic exchange matrix of finite type is of the form $B(c)$ for
some Coxeter element~$c$.
Thus Theorem~\ref{th:ca-general-1} provides a realization of any
cluster algebra of finite type with principal coefficients at an arbitrary
acyclic seed.

We next generalize Part (2) of
Theorem~\ref{th:ca-principal-minors-A-special}, by showing that
the set of all cluster variables in $\Acal(c)$ is of
the form $\{x_{\gamma;c}: \gamma \in \Pi(c)\}$ for some set of weights~$\Pi(c)$.
To describe~$\Pi(c)$, we need to develop some combinatorics related to the action on
weights by the cyclic group $\langle c \rangle$ generated by~$c$.
Interestingly, this combinatorics turns out to be very close to the one
developed in the classical paper \cite{kostant-1} by B.~Kostant, and
used recently in \cite{kirth} for (seemingly) completely
different purposes.

We consider the usual partial order on weights: $\gamma
\geq \delta$ if $\gamma - \delta$ is a nonnegative integer linear
combination of simple roots.
We also denote by $i \mapsto i^\star$ the involution on~$I$ given by
\begin{equation}
\label{eq:istar}
\omega_{i^\star} = - w_\circ (\omega_i),
\end{equation}
where, as usual, $w_\circ$ is the longest element of~$W$.

\begin{proposition}
\label{pr:hic}
For every~$i \in I$, there is a positive integer $h(i;c)$ such that
\begin{equation}
\label{eq:hic}
\omega_i > c \omega_i > c^{2} \omega_i > \cdots > c^{h(i;c)} \omega_i
= -\omega_{i^\star}.
\end{equation}
\end{proposition}

\begin{theorem}
\label{th:ca-general-2}
The set of cluster variables in $\Acal(c)$ is
$\{x_{\gamma;c}: \gamma \in \Pi(c)\}$, where
the set of weights $\Pi(c)$ is given by
\begin{equation}
\label{eq:Pi-c}
\Pi(c) = \{c^{m} \omega_i: i \in I, \, 0 \leq m \leq h(i;c)\}.
\end{equation}
\end{theorem}

As for Part (3) of Theorem~\ref{th:ca-principal-minors-A-special},
in the general case we give explicit formulas not for all exchange
relations but only for the \emph{primitive} ones, that is, those in which
one of the products of cluster variables on the right-hand side is
equal to~$1$.

\begin{theorem}
\label{th:ca-general-3}
The primitive exchange relations in $\Acal(c)$
are exactly the following:
\begin{equation}
\label{eq:exchange-primitive-special}
x_{-\omega_k;c}\ x_{\omega_k;c} =
y_{k;c} \prod_{i \prec_c k} x_{\omega_i;c}^{-a_{i,k}}
\prod_{k \prec_c i} x_{-\omega_i;c}^{-a_{i,k}} + 1;
\end{equation}
\begin{equation}
\label{eq:exchange-primitive-nonspecial}
x_{c^{m-1} \omega_k;c}\ x_{c^{m} \omega_k;c} =
\prod_{i \prec_c k} x_{c^{m} \omega_i;c}^{-a_{i,k}}
\prod_{k \prec_c i} x_{c^{m-1} \omega_i;c}^{-a_{i,k}} +
\prod_{j \in I} y_{j;c}^{[c^{m-1} \omega_k - c^{m} \omega_k :\ \alpha_j]},
\end{equation}
where $k \in I$, $1 \leq m \leq h(k;c)$, and
$[\gamma : \alpha_j]$ stands for the coefficient
of~$\alpha_j$ in the expansion of a weight~$\gamma$ in the
basis of simple roots.
\end{theorem}

The fact that the weights $c^{\,m-1} \omega_i$ and $c^{\,m} \omega_i$
appearing in the right hand side of \eqref{eq:exchange-primitive-nonspecial}
belong to $\Pi(c)$, is guaranteed by the following key property of
the numbers $h(i;c)$.

\begin{proposition}
\label{pr:hic-hjc}
If $i \prec_c j$ then we have
\begin{equation}
\label{eq:hic-hjc}
h(i;c) - h(j;c) = \begin{cases}
1 & \text{if $j^\star \prec_c i^\star$;} \\[.05in]
0 & \text{if $i^\star \prec_c j^\star$.}
\end{cases}
\end{equation}
\end{proposition}

Note that the numbers $h(i;c)$ are uniquely determined by
the relations \eqref{eq:hic-hjc} combined with the following property.

\begin{proposition}
\label{pr:hic+histarc}
For every $i \in I$, the sum
$h(i;c) + h(i^\star;c)$ is equal to the Coxeter number of the
connected component of~$I$ containing $i$ and $i^\star$.
\end{proposition}

In what follows, we identify the index set~$I$ with $[1,n] = \{1,
\dots, n\}$ so that the (fixed) Coxeter element~$c$ has the form
$c = s_1 \cdots s_n$.
Recall from \cite[(7.2),(7.3)]{ca4} that associated to every cluster variable~$z$
is the \emph{denominator vector} $\dd_{z;\xx} = (d_1, \dots, d_n) \in \ZZ^n$
with respect to a cluster~$\xx$.
Namely,~$z$ can be uniquely written in the form
\begin{equation}
\label{eq:Laurent-normal-form}
z = \frac{N(x_1, \dots, x_n)}{x_1^{d_1} \cdots x_n^{d_n}} \, ,
\end{equation}
where $N(x_1, \dots, x_n)$ is a polynomial not
divisible by any cluster variable~$x_i\in\xx$.
The following result gives an explicit formula
for the denominator vector of any
cluster variable~$z \in \Acal(c)$ with respect to the initial cluster.

\begin{theorem}
\label{th:denom-vector}
Let $\xx$ be the initial cluster in Theorem~\ref{th:ca-general-1},
and let $z = x_{\gamma;c}$ be a cluster variable not belonging to $\xx$.
Identifying $\ZZ^n$ with the root lattice by means of the
basis $\alpha_1, \dots, \alpha_n$ of simple roots,
the denominator vector $\dd_{z;\xx}$ gets identified with
$c^{-1}\gamma - \gamma$.
\end{theorem}

The following corollary agrees with \cite[Conjecture~7.4 (i),(ii)]{ca4}.

\begin{corollary}
\label{cor:denoms}
In the situation of Theorem~\ref{th:denom-vector}, all the
components~$d_i$ of the denominator vector $\dd_{z;\xx}$ are
nonnegative.
Furthermore, $d_i = 0$ if and only if~$z$ and the initial cluster
variable $x_{\omega_i;c}$ belong to the same cluster in $\Acal(c)$.
\end{corollary}

According to \cite[Corollary~6.3]{ca4}
(see also Section~\ref{sec:Cluster-algebras} below),
for any choice of coefficients,
expressing any cluster variable $z$
in terms of an initial cluster~$\xx$ only requires
knowing the \emph{$\gg$-vector} $\gg_{z;\xx}  \in
\ZZ^n$ and the \emph{$F$-polynomial} $F_{z;\xx} \in
\ZZ[t_1, \dots, t_n]$.
Our realization of cluster variables allows us to
obtain explicit expressions for their $\gg$-vectors and $F$-polynomials.
We start with $\gg$-vectors.

\begin{theorem}
\label{th:g-vector}
Let $\xx$ be the initial cluster in Theorem~\ref{th:ca-general-1},
and let $z = x_{\gamma;c}$ for some $\gamma \in \Pi(c)$.
Identifying $\ZZ^n$ with the weight lattice by means of the
basis $\omega_1, \dots, \omega_n$ of fundamental weights,
the $\gg$-vector $\gg_{z;\xx}$ gets identified with~$\gamma$.
\end{theorem}

\begin{remark}
\label{rem:cambrian-1}
An alternative description of these $\gg$-vectors was given in
\cite[Theorem~10.2]{reading-speyer}.
It is stated in different terms, and proved by a quite different
method, relying on \cite[Conjecture~7.12]{ca4}.
It is not difficult to check the equivalence of the two descriptions.
\end{remark}

To give a formula for the $F$-polynomials, recall that, for each
$i \in [1,n]$, there are one-parameter root subgroups in~$G$ given by
\begin{equation}
\label{eq:xi-ibar}
x_i (t) = \varphi_i \mat{1}{t}{0}{1}, \quad x_{\bar i} (t) = \varphi_i \mat{1}{0}{t}{1},
\end{equation}
where $\varphi_i: SL_2 \to G$ denotes
the canonical embedding corresponding to the simple
root~$\alpha_i$.

\begin{theorem}
\label{th:F-poly}
Under the assumptions of Theorem~\ref{th:g-vector},
the $F$-polynomial $F_{z;\xx}$ is given by
\begin{equation}
\label{eq:F-poly}
F_{z;\xx}(t_1, \dots, t_n) = \Delta_{\gamma, \gamma}
(x_{\bar 1}(1) \cdots x_{\bar n}(1) x_{n}(t_n) \cdots x_{1}(t_1)).
\end{equation}
\end{theorem}

The following corollary agrees with \cite[Conjecture~5.4]{ca4}.

\begin{corollary}
\label{cor:ct-1}
Every $F$-polynomial in Theorem~\ref{th:F-poly} has constant term~$1$.
\end{corollary}

\begin{example}
\label{ex:c-An-special}
Let $c$ be as in Theorem~\ref{th:ca-principal-minors-A-special}
(for $G$ of type~$A_n$).
Then we have $h(k;c) = n+1-k$ for $k = 1, \cdots, n$.
By Theorem~\ref{th:g-vector}, the cluster variable $x_{c^m
\omega_k;c} = x_{[m+1,m+k]}$ has the $\gg$-vector
$c^m \omega_k = \omega_{m+k} - \omega_m$ (with the convention
$\omega_0 = \omega_{n+1} = 0$).
By Theorem~\ref{th:denom-vector}, for $m \geq 1$, the denominator
vector of $x_{c^m \omega_k;c}$ is equal to
$$c^{m-1} \omega_k - c^m \omega_k = \alpha_m + \alpha_{m+1} +
\cdots + \alpha_{m+k-1}.$$
By computing the determinant in Theorem~\ref{th:F-poly}, we conclude that $x_{c^m
\omega_k;c}$ has the $F$-polynomial
$$F(t_1, \dots, t_n) = 1 + t_m + t_m t_{m+1} + \cdots +
t_m t_{m+1} \cdots t_{m+k-1}$$
(with the convention $t_0 = 0$).
\end{example}

The rest of the paper is devoted to the proofs of the above results.
Most of our proofs rely on the following essentially
combinatorial fact: every Coxeter element can be obtained from
any other by a sequence of operations sending $c = s_{i_1} \cdots s_{i_n}$ to
$\tilde c = s_{i_1} c s_{i_1}$.
Thus to prove some statement for an arbitrary Coxeter element~$c$, it
is enough to check that it holds for some special choice of~$c$,
and that it is preserved by the above operations.
As this special choice, we take a \emph{bipartite} Coxeter
element~$t$, for which most of the above results were essentially
established in earlier papers (see Section~\ref{sec:cce} for more details).

In Section~\ref{sec:cce}, after a short reminder on root systems, we develop the
combinatorics of the action of a Coxeter element on roots and
fundamental weights, in particular proving
Propositions~\ref{pr:hic},~\ref{pr:hic-hjc}, and~\ref{pr:hic+histarc}.
In Section~\ref{sec:primitive-exchange-relations}, we recall the
necessary background on (reduced) double Bruhat cells and
generalized minors, and prove the relations \eqref{eq:exchange-primitive-special}
and \eqref{eq:exchange-primitive-nonspecial} from
Theorem~\ref{th:ca-general-3}.
Both sections~\ref{sec:cce} and~\ref{sec:primitive-exchange-relations}
are totally independent of the theory of cluster algebras, which
are not even mentioned there.

In Section~\ref{sec:Cluster-algebras} we recall
basic definitions and facts about cluster algebras, following mainly
\cite{ca4}, then prove Theorem~\ref{th:ca-general-1}.
Note that some of the preliminary results given there
(Propositions~\ref{pr:principal-embedding}, \ref{pr:trop-semifield-homs}
and \ref{pr:coef-specialization}) have not appeared in this form
before, although they are easy consequences of the known results.

In Section~\ref{sec:univcoeff} we prove
Theorems~\ref{th:ca-general-2} and~\ref{th:ca-general-3}.
In the course of the proof we obtain some results of independent interest.
Proposition~\ref{pr:comp-degree-Pic} transfers the compatibility degree function
introduced in \cite[Section~3]{yga} to each set~$\Pi(c)$.
This leads to a combinatorial description of all the clusters in~$\Acal(c)$
given in Corollary~\ref{cor:c-clusters}.
Another important result is Proposition~\ref{pr:universal-ca-c},
which provides, for every choice of a Coxeter element~$c$, a
realization of the cluster algebra with \emph{universal
coefficients} introduced in \cite[Section~12]{ca4}, with respect to the initial seed
with the exchange matrix $B(c)$.
As a corollary of the technique developed in this section, we
describe explicit isomorphisms between various algebras
$\Acal(c)$ (Corollary~\ref{cor:c-tildec-isom-algebraic})
as well as corresponding geometric isomorphisms between
various reduced double cells $L^{c,c^{-1}}$
(Remark~\ref{rem:c-tildec-isom-geometric}).

In Section~\ref{sec:proofs-g-vector-F-poly}, we prove
Theorems~\ref{th:denom-vector}, \ref{th:g-vector} and \ref{th:F-poly},
as well as Corollaries~\ref{cor:denoms} and \ref{cor:ct-1}.
As explained in Remark~\ref{rem:cambrian-2}, these results allow
us to give an alternative combinatorial description of the
$c$-Cambrian fans studied in \cite{reading-speyer}.

In the final
Section~\ref{sec:proof-of-ca-principal-minors-A-special}, we
illustrate the results of the paper by an example dealing with
type~$A_n$ and the Coxeter element~$c = s_1 \cdots s_n$ in the
standard numbering of simple roots.
In particular, we prove Theorem~\ref{th:ca-principal-minors-A-special}
and present an example making more explicit the construction given in
Remark~\ref{rem:c-tildec-isom-geometric}.

\section{Proofs of Propositions~\ref{pr:hic},~\ref{pr:hic-hjc}, and~\ref{pr:hic+histarc}}
\label{sec:cce}

We start by laying out the basic terminology and notation related
to root systems, to be used throughout the paper;
some of it has already appeared in the introduction.
In what follows, $A = (a_{i,j})_{i,j \in I}$ is an indecomposable
$n\times n$ Cartan matrix, i.e., one of the matrices $A_n, B_n, \dots, G_2$ in the
Cartan-Killing classification (the general case in
Propositions~\ref{pr:hic},~\ref{pr:hic-hjc}, and~\ref{pr:hic+histarc}
easily reduces to the case of~$A$ indecomposable).

Let $\Phi$ be the corresponding rank $n$
root system with the set of simple roots $\{\alpha_i: i \in I\}$.
Let $W$ be the Weyl group of $\Phi$, i.e., the group of linear
transformations of the root space, generated by the \emph{simple
reflections} $s_i \ (i \in I)$ whose action on simple roots is
given by
\begin{equation}
\label{eq:si-alphaj}
s_i (\alpha_j) = \alpha_j - a_{i,j} \alpha_i.
\end{equation}
Let $\{\omega_i: i \in I\}$ be the set of fundamental weights
related to simple roots via
\begin{equation}
\label{eq:alpha-omega}
\alpha_j = \sum_{i \in I} a_{i,j} \omega_i.
\end{equation}
The action of  simple reflections on the fundamental weights is given by
\begin{equation}
\label{eq:si-omegaj}
s_i \omega_j = \begin{cases}
\omega_i - \alpha_i & \text{if $j = i$;} \\[.05in]
\omega_j & \text{if $j \neq i$.}
\end{cases}
\end{equation}

For future use, here are a couple of useful lemmas.

\begin{lemma}
\label{lem:two-identities}
Suppose $I = \{1, \dots, n\}$ and $c = s_1 \cdots s_n$.
For $i = 1, \dots, n$, let
\begin{equation}
\label{eq:beta-i}
\beta_i = s_1 \cdots s_{i-1} \alpha_i,
\end{equation}
so that $\{\beta_1, \dots, \beta_n\}$ is the set of positive
roots~$\beta$ such that $c^{-1} \beta$ is negative.
Then we have
\begin{equation}
\label{eq:omega-c-omega}
\omega_i - c \omega_i = \beta_i,
\end{equation}
and
\begin{equation}
\label{eq:beta-telescoping}
\beta_j + \sum_{i =1}^{j-1} a_{i,j} \beta_i = \alpha_j.
\end{equation}
\end{lemma}

\begin{proof}
The identity (\ref{eq:omega-c-omega}) is an immediate consequence of
\eqref{eq:si-omegaj}:
$$\beta_i
= s_1\cdots s_{i-1}(\omega_i - s_i\omega_i) = \omega_i - s_1\cdots s_i\omega_i
= \omega_i - c\omega_i.$$
As for (\ref{eq:beta-telescoping}), we have
$$\beta_j + \sum_{i =1}^{j-1} a_{i,j} \beta_i =
s_1\cdots s_{j-1}\alpha_j + \sum_{i =1}^{j-1} s_1\cdots
s_{i-1}(\alpha_j - s_i \alpha_j) = \alpha_j,$$
by telescoping.
\end{proof}

\begin{lemma} \cite[Proposition~VI.1.33]{bourbaki}
\label{lem:c-reps-in-PHI}
In the situation of Lemma~\ref{lem:two-identities}, let
$\langle c \rangle$ denote the cyclic subgroup of~$W$ generated by~$c$.
Then every $\langle c \rangle$-orbit in~$\Phi$ contains exactly
one positive root~$\beta$ such that $c^{-1} \beta$ is negative
(i.e., one of the roots $\beta_1, \dots, \beta_n$) and exactly
one positive root~$\alpha$ such that $c \alpha$ is negative
(i.e., one of the roots $- c^{-1} \beta_1, \dots, - c^{-1} \beta_n$).
\end{lemma}

A \emph {reduced word} for $w\in W$ is a sequence of indices
$\ii=(i_1,\dots,i_\ell)$ of shortest possible length $\ell = \ell(w)$
such that $w=s_{i_1} \cdots s_{i_\ell}$.
The group $W$ possesses a unique longest element denoted by~$w_\circ$.

The \emph{Coxeter graph} associated to $\Phi$ has the index set $I$
as the set of vertices, with $i$ and $j$ joined by an
edge whenever $a_{i,j} a_{j,i} > 0$.
Since we assume that $A$ is indecomposable,
the root system $\Phi$ is irreducible,
and the Coxeter graph $I$ is a tree
(one of the familiar Coxeter-Dynkin diagrams $A_n, D_n, E_6, E_7, E_8$).

Recall that a \emph{Coxeter element} is the product of all simple
reflections $s_i$ (for $i \in I$) taken in an arbitrary order.
Thus, we have $\ell(c) = n$.
It is well-known that all Coxeter elements are conjugate to each
other in~$W$; in particular, they have the same order~$h$ called the
\emph{Coxeter number} of~$W$.

Note that Coxeter elements are in a natural bijection with
orientations of the Coxeter graph: we orient an edge $j \to i$ if
$i \prec_c j$, i.e., if $i$ precedes $j$ in some (equivalently,
any) reduced word for~$c$.
We introduce the following elementary operation on Coxeter
elements:
\begin{align}
\label{eq:cyclical-move}
&\text{replace~$c$ with $\tilde c$ if $c = s_i c'$ and $\tilde c =
c' s_i$ for}\\
\nonumber
&\text{some $i \in I$ and $c' \in W$ with $\ell(c') = n-1$.}
\end{align}
We will use the following well-known fact:
\begin{align}
\label{eq:cyclical-moves-transitive}
&\text{any Coxeter element can be reached from}\\
\nonumber
&\text{any other one by a
sequence of moves \eqref{eq:cyclical-move}.}
\end{align}
(Passing from Coxeter elements to corresponding orientations of
the Coxeter graph, \eqref{eq:cyclical-moves-transitive} translates
into the statement that any orientation of a tree can be
obtained from any other orientation by a repeated application of
the following operation: reversing all arrows at some sink; this
is easily proved by induction on the size of a tree.)

Our proofs of the needed properties of Coxeter elements will use
the same strategy: show that the property in question is preserved
under any operation \eqref{eq:cyclical-move}, and prove that it
holds for some particular choice of a Coxeter element.
As this particular choice, we use the \emph{bipartite} Coxeter
element defined as follows.

Since the Coxeter graph is a tree, it is bipartite, so the set of
vertices $I$ is a disjoint union of two parts
$I_+$ and $I_-$ such that every edge joins two vertices from
different parts. Note that $I_+$ and $I_-$ are determined uniquely up to renaming.
We write $\varepsilon (i) = \varepsilon$ for $i \in I_\varepsilon$.

Now we define the bipartite Coxeter element $t$ by setting
\begin{equation}
\label{eq:coxeter-bipartite}
t = t_+ t_-,
\end{equation}
where
\begin{equation}
\label{eq:coxeter-t-pm}
t_\pm = \prod_{\varepsilon(i) = \pm 1} s_i\,.
\end{equation}
Note that the order of factors in (\ref{eq:coxeter-t-pm}) does not matter because
$s_i$ and $s_j$ commute whenever $\varepsilon (i) = \varepsilon (j)$.
Let $\ii_-$ and $\ii_+$ be some reduced words for the elements $t_-$ and~$t_+\,$, respectively.

\begin{lemma}
\label{lem:lusztig-rw}
\cite[Exercise~V.6.2]{bourbaki}
We have $\underbrace{t_+ t_- t_+ \cdots t_\pm
  t_\mp}_{\text{$h$ factors}}=w_\circ\,$.
Moreover,
the word
\begin{equation}
\label{eq:special-rw}
\ii_\circ \stackrel{\rm def}{=}
\underbrace{\ii_+ \ii_- \ii_+ \cdots \ii_\pm
  \ii_\mp}_h
\end{equation}
(concatenation of $h$ segments)
is a reduced word for~$w_\circ\,$.
\end{lemma}

Regarding Lemma~\ref{lem:lusztig-rw}, recall from the tables in
\cite{bourbaki} that $h$ is even for all types except $A_n$ with $n$
even; in the exceptional case of type $A_{2e}$, we have $h = 2e +1$.
If $h$ is even, then Lemma~\ref{lem:lusztig-rw} says that $t^{h/2} = \wnot$.

Now everything is in place for the proofs of
Propositions~\ref{pr:hic},~\ref{pr:hic-hjc},
and~\ref{pr:hic+histarc}.
We start by observing that Proposition~\ref{pr:hic} is a consequence
of the following seemingly weaker statement:
\begin{align}
\label{eq:omegai-w0omegaistar}
&\text{for every Coxeter element~$c$ and $i \in I$, the}\\
\nonumber
&\text{weight $w_\circ \omega_i$ belongs to the $\langle c \rangle$-orbit of $\omega_i$.}
\end{align}
Indeed, assuming \eqref{eq:omegai-w0omegaistar}, we can define
$h(i;c)$ as the smallest positive integer~$n$ such that
$c^m \omega_i = w_\circ \omega_i$.
To prove the inequalities in \eqref{eq:hic}, note that, in view of
\eqref{eq:omega-c-omega}, we have
$c^m \omega_i - c^{m+1} \omega_i = c^m \beta_i$
for $m \in \ZZ$, and so the differences $c^m \omega_i - c^{m+1} \omega_i$
form the $\langle c \rangle$-orbit of $\beta_i$ in $\Phi$.
Now recall Lemma~\ref{lem:c-reps-in-PHI}.
Since $\omega_i$ is the maximal element of its $W$-orbit, while
$w_\circ \omega_i$ is the minimal element of the same orbit, we
see that $\beta = \omega_i - c \omega_i$ (resp.~$\alpha = c^{-1} w_\circ \omega_i -
w_\circ \omega_i$) is the unique positive root in the $\langle c \rangle$-orbit of $\beta_i$
such that $c^{-1} \beta$ (resp.~$c \alpha$) is negative.
It follows that all the roots $\gamma = c^m \omega_i - c^{m+1}
\omega_i$ with $0 \leq m < h(i;c)$ are positive, as required.

As an easy consequence of Lemma~\ref{lem:lusztig-rw}, the
statement \eqref{eq:omegai-w0omegaistar} (hence Proposition~\ref{pr:hic})
holds for the bipartite Coxeter element~$t$; furthermore, we have
\begin{equation}
\label{eq:h(i;t)}
h(i;t) = \begin{cases}
\lfloor \frac{h}{2} \rfloor  & \text{for $\varepsilon(i) = -1$\ ;} \\[.05in]
\lceil \frac{h}{2} \rceil  & \text{for $\varepsilon(i) = +1$\ .}
\end{cases}
\end{equation}
The fact that~$t$ satisfies Propositions~\ref{pr:hic-hjc}
and~\ref{pr:hic+histarc}, follows at once from \eqref{eq:h(i;t)}
together with an observation that, for every $i \in I$, we have
\begin{equation}
\label{eq:istar-parts}
\varepsilon (i^\star) = \begin{cases}
\varepsilon(i)  & \text{if $h$ is even;} \\[.05in]
-\varepsilon(i)  & \text{if $h$ is odd.}
\end{cases}
\end{equation}

In view of \eqref{eq:cyclical-moves-transitive}, to finish the proofs of
Propositions~\ref{pr:hic},~\ref{pr:hic-hjc},
and~\ref{pr:hic+histarc}, it suffices to assume that they hold
for some Coxeter element~$c$ and show that the same is true for
the element~$\tilde c$ obtained from~$c$ via \eqref{eq:cyclical-move}.
Without loss of generality, we assume that $I = \{1, \dots, n\}$,
$c = s_1 s_2 \cdots s_n$, and $\tilde c = s_2 \cdots s_n s_1$.

The key statement to prove is the following: $\tilde c$ satisfies
\eqref{eq:omegai-w0omegaistar} (hence Proposition~\ref{pr:hic}),
and we have
\begin{equation}
\label{eq:h(i;tildec)}
h(i;\tilde{c}) = \begin{cases}
h(i;c) - 1 & \text{if $i = 1$ and $i^\star \neq 1$;} \\[.05in]
h(i;c) + 1 & \text{if $i \neq 1$ and $i^\star = 1$;} \\[.05in]
h(i;c) & \text{otherwise.}
\end{cases}
\end{equation}

We split the proof of \eqref{eq:h(i;tildec)} into four cases:

\smallskip

{\bf Case~1.} Let $i \neq 1$ and $i^\star \neq 1$.
Since $\tilde c^m = s_1 c^m s_1$ for all $m > 0$, using
\eqref{eq:si-omegaj}, we obtain
$$\tilde c^m \omega_i = - \omega_{i^\star} \Longleftrightarrow
s_1 c^m s_1 \omega_i = - \omega_{i^\star} \Longleftrightarrow
c^m  \omega_i = - \omega_{i^\star},$$
implying that $h(i;\tilde{c}) = h(i;c)$.

\smallskip

{\bf Case~2.} Let $i = i^\star = 1$.
Writing $\tilde c^m = s_2 \cdots s_n c^m s_n \cdots s_2$ and using
\eqref{eq:si-omegaj}, we obtain
$$\tilde c^m \omega_1 = - \omega_{1}  \Longleftrightarrow
c^m  \omega_1 = - \omega_{1},$$
again implying that $h(1;\tilde{c}) = h(1;c)$.

\smallskip

{\bf Case~3.} Let $i = 1$ and $i^\star \neq 1$.
Writing $\tilde c^m = s_1 c^{m+1} s_n \cdots s_2$ and using
\eqref{eq:si-omegaj}, we obtain
$$\tilde c^m \omega_1 = - \omega_{i^\star} \Longleftrightarrow
c^{m+1}  \omega_1 = - \omega_{i^\star},$$
implying that $h(i;\tilde{c}) = h(i;c) - 1$.

\smallskip

{\bf Case~4.} Finally, let $i \neq 1$ and $i^\star = 1$.
Writing $\tilde c^m = s_2 \cdots s_n c^{m-1} s_1$ and using
\eqref{eq:si-omegaj}, we obtain
$$\tilde c^m \omega_i = - \omega_{1} \Longleftrightarrow
c^{m-1}  \omega_i = - \omega_{1},$$
implying that $h(i;\tilde{c}) = h(i;c) + 1$.

\smallskip

This concludes the proof of Proposition~\ref{pr:hic} and the
relation \eqref{eq:h(i;tildec)}.

\medskip

The proofs of Propositions~\ref{pr:hic-hjc}
and~\ref{pr:hic+histarc} now become purely combinatorial exercises.
Namely, let again $c = s_1 s_2 \cdots s_n$ and $\tilde c = s_2 \cdots s_n s_1$.
To prove Proposition~\ref{pr:hic+histarc}, it suffices to show
that $h(i;\tilde{c}) + h(i^\star;\tilde{c}) = h(i;c) +
h(i^\star;c)$ for all~$i$, which is immediate from
\eqref{eq:h(i;tildec)}.

To prove Proposition~\ref{pr:hic-hjc}, it suffices to show the
following: if the numbers $h(i;c)$ satisfy \eqref{eq:hic-hjc}, and
the $h(i;\tilde c)$ are given by \eqref{eq:h(i;tildec)}, then the
$h(i;\tilde c)$ also satisfy \eqref{eq:hic-hjc} with $c$ replaced by $\tilde c$.
There are several cases to consider.

\smallskip

{\bf Case~1.} Suppose that $i \prec_{\tilde{c}} j$, and $i^\star \prec_{\tilde{c}} j^\star$.
Then we have $i \neq 1$ and $i^\star \neq 1$, hence $h(i;\tilde c) = h(i;c)$.
This case breaks into four subcases according to which of the two
indices $j$ and $j^\star$ are equal to~$1$.
In each of these subcases, the desired equality $h(j;\tilde c) = h(i;\tilde c)$
is seen by a direct inspection.
For instance, if $j = 1 \neq j^\star$, we have $h(j; \tilde c) =
h(j;c)-1$ by \eqref{eq:h(i;tildec)}; on the other hand, in this
case we have $j \prec_c i$, and $i^\star \prec_c j^\star$, hence
$h(j;c)-1= h(i;c)$ by \eqref{eq:hic-hjc}.

\smallskip

{\bf Case~2.} Now suppose that $i \prec_{\tilde{c}} j$, and
$j^\star \prec_{\tilde{c}} i^\star$.
Then we have $i \neq 1$ and $j^\star \neq 1$.
This case also breaks into four subcases according to which of the two
indices $j$ and $i^\star$ are equal to~$1$.
Again, in each of these subcases, the desired equality $h(j;\tilde c) = h(i;\tilde c)-1$
is seen by a direct inspection.
For instance, if $j = 1 \neq i^\star$, we have $h(i;\tilde c) = h(i;c)$, and
$h(j; \tilde c) = h(j;c)-1$ by \eqref{eq:h(i;tildec)}; on the other hand, in this
case we have $j \prec_c i$, and $j^\star \prec_c i^\star$, hence
$h(j;c)= h(i;c)$ by \eqref{eq:hic-hjc}.

\medskip

This concludes the proofs of Propositions~\ref{pr:hic},~\ref{pr:hic-hjc}, and~\ref{pr:hic+histarc}.

\section{Proofs of relations \eqref{eq:exchange-primitive-special} and \eqref{eq:exchange-primitive-nonspecial}}
\label{sec:primitive-exchange-relations}

For the convenience of the reader, we start by briefly recalling
necessary facts about generalized minors and reduced double cells;
more details can be found in \cite{fz-double,bz01}.


Let $\lg$ be a complex semisimple Lie algebra of rank~$n$ with Chevalley generators
$f_i, \alpha_i^\vee$, and $e_i$ for $i \in I$, where $I$ is an $n$-element index set, which will be often
identified with $\{1, \dots, n\}$.
The elements $\alpha_i^\vee$ are \emph{simple coroots} of $\lg$;
they form a basis of a Cartan subalgebra $\hh$ of $\lg$.
The \emph{simple roots} $\alpha_i \,\, (i \in I)$ form a basis in the dual space
$\hh^*$ such that $[h, e_i] = \alpha_i (h) e_i$, and $[h,f_i] = - \alpha_i (h) f_i$
for any $h \in \hh$ and $i \in I$.
The structure of $\gg$ is uniquely determined by the
\emph{Cartan matrix} $A = (a_{i,j})$ given by $a_{ij} = \alpha_j (\alpha_i^\vee)$.

Let $G$ be the simply connected complex semisimple Lie group with the Lie algebra $\lg$.
For every $i\in I$, let $\varphi_i: SL_2 \to G$ denote the canonical embedding
corresponding to the simple root~$\alpha_i\,$.
We use the notation
\begin{equation}
\label{eq:xi-xibar}
x_i (t) = \varphi_i \mat{1}{t}{0}{1} = \exp \ (t e_i), \quad
x_{\bar i} (t) = \varphi_i \mat{1}{0}{t}{1} = \exp \ (t f_i) \ .
\end{equation}
We also set
$$t^{\alpha_i^\vee} = \varphi_i \mat{t}{0}{0}{t^{-1}} \in H$$
for any $i \in I$ and any $t \neq 0$.
Let $N$ (resp. $N_-$) be the maximal unipotent subgroup of $G$ generated by all $x_i (t)$ (resp. $x_{\bar i} (t)$).
Let $H$ be the maximal torus in $G$ with the Lie algebra $\hh$.
Let $B = HN$ and $B_- = HN_-$ be the corresponding two opposite Borel subgroups.

The Weyl group $W$ associated to~$A$ is naturally identified with ${\rm Norm}_G (H)/H$;
this identification sends each simple reflection
$s_i$ to the coset $\overline {s_i} H$, where the representative $\overline {s_i} \in {\rm Norm}_G (H)$ is defined by
\[
\overline {s_i} = \varphi_i \mat{0}{-1}{1}{0} \, .
\]
The elements $\overline {s_i}$ satisfy the braid relations in~$W$;
thus the representative $\overline w$ can be unambiguously defined for any
$w \in W$ by requiring that $\overline {uv} = \overline {u} \cdot \overline {v}$
whenever $\l (uv) = \l (u) + \l (v)$.

With some abuse of notation, we identify the weight lattice $P$ in $\hh^*$
with the group of rational multiplicative characters of~$H$,
here written in the exponential notation:
a weight $\gamma\in P$ acts by $a \mapsto a^\gamma$.
Under this identification, the fundamental weights $\omega_1, \ldots, \omega_n$ act in $H$ by
$(t^{\alpha_j^\vee})^{\omega_i} = t^{\delta_{ij}}$.

Recall that the set $G_0 = N_-HN$ of elements $x\in G$ that
have Gaussian decomposition is open and dense in~$G$; this (unique) decomposition of
$x \in N_-HN$ is written as $x = [x]_- [x]_0 [x]_+ \,$.


We now recall some basic properties of \emph{generalized minors} introduced in \cite{fz-double}.
For $u,v \in W$ and $i \in I$, the \emph{generalized minor}
$\Delta_{u \omega_i, v \omega_i}$
is the regular function on $G$ whose restriction to the open set
${\overline {u}} G_0 {\overline {v}}^{-1}$ is given by
\begin{equation}
\label{eq:Delta-general}
\Delta_{u \omega_i, v \omega_i} (x) =
(\left[{\overline {u}}^{\ -1}
   x \overline v\right]_0)^{\omega_i} \ .
\end{equation}
As shown in \cite{fz-double}, $\Delta_{u \omega_i, v \omega_i}$ depends on
the weights $u \omega_i$ and $v \omega_i$ alone, not on the particular
choice of $u$ and~$v$.
In the special case $G=SL_{n+1}\,$, the generalized minors are nothing but
the ordinary minors of a matrix.

Generalized minors have the following properties
(see \cite[(2.14), (2.25)]{fz-double}):
\begin{equation}
\label{eq:minor-weight}
\Delta_{\gamma, \delta} (a_1 x a_2) = a_1^\gamma a_2^\delta
\Delta_{\gamma, \delta} (x) \quad (a_1, a_2  \in H; \ x \in G)\ ;
\end{equation}
\begin{equation}
\label{eq:minor-transpose}
\Delta_{\gamma,\delta} (x) =  \Delta_{\delta, \gamma} (x^T) \ ,
\end{equation}
where $x \mapsto x^T$ is the ``transpose" involutive antiautomorphism of~$G$
acting on generators by
\begin{equation}
\label{eq:T}
a^T = a \quad (a \in H) \ , \quad x_i (t)^T = x_{\bar i} (t) \ ,
\quad x_{\bar i} (t)^T = x_i (t) \ .
\end{equation}
We will also use the involutive antiautomorphism $\iota$ of $G$
introduced in \cite[(2.2)]{fz-double}; it is defined by
\begin{equation}
\label{eq:iota}
a^\iota = a^{-1} \quad (a \in H) \ , \quad x_i (t)^\iota = x_i (t) \ ,
\quad x_{\bar i} (t)^\iota = x_{\bar i} (t) \ .
\end{equation}
By \cite[(2.25)]{fz-double}, we have
\begin{equation}
\label{eq:minors-tau}
\Delta_{\gamma, \delta} (x) = \Delta_{-\delta, -\gamma} (x^\iota)
\end{equation}
for any generalized minor $\Delta_{\gamma, \delta}$, and any $x \in G$.

The following important identity was obtained in \cite[Theorem~1.17]{fz-double}.

\begin{proposition}
\label{pro:minors-Dodgson}
Suppose $u,v \in W$ and $k \in I$
are such that $us_k \omega_k < u \omega_k$ and $vs_k \omega_k < v \omega_k$.
Then
\begin{eqnarray}
\begin{array}{l}
\label{eq:minors-Dodgson}
\Delta_{u \omega_k, v \omega_k} \Delta_{us_k \omega_k, v s_k \omega_k}
= \Delta_{us_k \omega_k, v \omega_k} \Delta_{u \omega_k, v s_k \omega_k}
+ \prod_{i \neq k} \Delta_{u \omega_i, v \omega_i}^{- a_{i,k}} \ .
\end{array}
\end{eqnarray}
\end{proposition}

The group $G$ has two \emph{Bruhat decompositions}, with respect to opposite Borel subgroups $B$ and $B_-\,$:
$$G = \bigcup_{u \in W} B u B = \bigcup_{v \in W} B_- v B_-  \ . $$
The \emph{double Bruhat cells}~$G^{u,v}$ are defined by $G^{u,v} = B u B  \cap B_- v B_- \,$.
These varieties were introduced and studied in \cite{fz-double}.
Each double Bruhat cell can be defined inside~$G$ by a collection of vanishing/non-vanishing conditions
of the form $\Delta(x)=0$ and \hbox{$\Delta(x)\neq0$,}
where $\Delta$ is a generalized minor.
The following description can be found in
\cite[Proposition~2.8]{ca3}.

\begin{proposition}
\label{pr:double-cell-by-eqs}
A double Bruhat cell $G^{u,v}$ is given inside~$G$ by the following
conditions, for all $i \in I$:
\begin{align}
& \Delta_{u'\omega_i,\omega_i}=0 \quad \text{whenever $u'\omega_i\not\leq
  u\omega_i$ in the Bruhat order,}\\
& \Delta_{\omega_i,v'\omega_i}=0 \quad \text{whenever $v'\omega_i\not\leq
  v^{-1}\omega_i$ in the Bruhat order,}\\
\label{eq:double-cell-ineq}
&\Delta_{u\omega_i,\omega_i}\neq 0, \quad
  \Delta_{\omega_i,v^{-1}\omega_i}\neq 0\,.
\end{align}
\end{proposition}

In this paper we concentrate on the following subset $L^{u,v} \subset
G^{u,v}$ introduced in \cite{bz01} and called a \emph{reduced double Bruhat cell}:
\begin{equation}
\label{eq:reduced cell}
L^{u,v} = N \overline u N  \cap B_- v B_- \ .
\end{equation}
The equations defining $L^{u,v}$ inside $G^{u,v}$ look as follows.

\begin{proposition}\cite[Proposition 4.3]{bz01}
\label{pr:Luv-equations}
An element $x \in G^{u,v}$ belongs to $L^{u,v}$ if and only if
$\Delta_{u \omega_i, \omega_i} (x) = 1$ for $i \in  [1,n]$.
\end{proposition}

The variety $L^{u,v}$ (resp. $G^{u,v}$) is biregularly isomorphic to a Zariski
open subset of an affine space of dimension $\l(u)+\l(v)$ (resp. $\l(u)+\l(v)+n$).
Local coordinates in $L^{u,v}$ and $G^{u,v}$ were constructed in \cite{fz-double,bz01}.
To simplify the notation, we will describe them only in the
following special case.

In the rest of the section we assume that the index set $I$ is
$\{1, \dots, n\}$, and work with the fixed Coxeter element
$c = s_1 \cdots s_n$, and with the varieties $L^{c, c^{-1}}$ and
$G^{c, c^{-1}}$.

In accordance with \cite[Theorem~1.2]{fz-double}, a generic element
$x \in G^{c, c^{-1}}$ can be uniquely factored as
\begin{equation}
\label{eq:Gc-c-1-factor}
x = x_{\bar 1}(u_1) \cdots x_{\bar n}(u_n)\, a\,
x_n(t_n) \cdots x_1(t_1)
\end{equation}
for some $a \in H$ and $u_i, t_i \in \CC^*$ (see
\eqref{eq:xi-xibar}).

In the case of the reduced double cell $L^{c, c^{-1}}$, the
factorization \eqref{eq:Gc-c-1-factor} can be modified as follows.
For any nonzero $u \in \CC$ and any $i = 1, \dots, n$, denote
\begin{equation}
\label{eq:xnegative}
x_{- i} (u) = x_{\bar i} (u) u^{- \alpha_i^\vee}  = \varphi_i \mat{u^{-1}}{0}{1}{u} \, .
\end{equation}
In accordance with  \cite[Proposition~4.5]{bz01}, a generic element
$x \in L^{c, c^{-1}}$ can be uniquely factored as
\begin{equation}
\label{eq:Lc-c-1-factor}
x = x_{-1}(u_1) \cdots x_{-n}(u_n) x_n(t_n) \cdots x_1(t_1),
\end{equation}
where $u_i, t_i \in \CC^*$.

After this preparation, let us turn to the proof of
\eqref{eq:exchange-primitive-nonspecial}.
Let $k = 1, \dots, n$ and $1 \leq m \leq h(k;c)$.
Specializing \eqref{eq:minors-Dodgson} for
$u = v = c^{m-1} s_1\cdots s_{k-1}$ (and taking into account \eqref{eq:si-omegaj}), we get
\begin{equation}
\label{eq:ns-minors-Dodgson}
\begin{split}
\Delta_{c^{m-1} \omega_k, c^{m-1} \omega_k} \Delta_{c^{m} \omega_k, c^{m} \omega_k} &=
\Delta_{c^{m} \omega_k, c^{m-1} \omega_k} \Delta_{c^{m-1} \omega_k, c^{m} \omega_k} \\
&+ \prod_{i < k} \Delta_{c^{m} \omega_i, c^{m} \omega_i}^{- a_{i,k}} \prod_{i > k}
\Delta_{c^{m-1} \omega_i, c^{m-1} \omega_i}^{- a_{i,k}}\ .
\end{split}
\end{equation}
Remembering \eqref{eq:x-gamma-c} and \eqref{eq:y-j-c}, and
comparing \eqref{eq:ns-minors-Dodgson} with the desired equality
\eqref{eq:exchange-primitive-nonspecial}, we see that it remains
to prove the following lemma.

\begin{lemma}
\label{lem:coeffs-nonspecial}
For every $k = 1, \dots, n$ and $1 \leq m \leq h(k;c)$, the
restriction to $L^{c,c^{-1}}$ of the function
$\Delta_{c^{m} \omega_k, c^{m-1} \omega_k} \Delta_{c^{m-1} \omega_k, c^{m} \omega_k}$
is equal to that of the product
\begin{equation}
\label{eq:product-coeffs}
\prod_{j=1}^n (\Delta_{\omega_j, c \omega_j} \prod_{i < j}
\Delta_{\omega_i, c \omega_i}^{a_{i,j}})^{[c^{m-1} \omega_k - c^{m} \omega_k :
\alpha_j]}.
\end{equation}
\end{lemma}

\begin{proof}
First of all, in view of \eqref{eq:omega-c-omega}, we have
\begin{equation}
\label{eq:cm-beta-k}
c^{m-1} \omega_k - c^{m} \omega_k = c^{m-1} \beta_k.
\end{equation}
Using \eqref{eq:beta-telescoping}, we can rewrite the product in
\eqref{eq:product-coeffs} as
\begin{equation}
\label{eq:product-coeffs-simplified}
\prod_{j=1}^n (\Delta_{\omega_j, c \omega_j} \prod_{i < j}
\Delta_{\omega_i, c \omega_i}^{a_{i,j}})^{[c^{m-1} \omega_k - c^{m} \omega_k :
\alpha_j]} = \prod_{i=1}^n \Delta_{\omega_i, c \omega_i}^{[c^{m-1} \beta_k : \beta_i]},
\end{equation}
where the exponent $[c^{m-1} \beta_k : \beta_i]$ stands for the
coefficient of $\beta_i$ in the expansion of $c^{m-1} \beta_k$ in the
basis $\{\beta_1, \dots, \beta_n\}$.

Now let us deal with the functions $\Delta_{c^{m} \omega_k, c^{m-1} \omega_k}$
and $\Delta_{c^{m-1} \omega_k, c^{m} \omega_k}$.
Let~$x$ be a generic element of $G^{c, c^{-1}}$ expressed as in \eqref{eq:Gc-c-1-factor}.
We claim that
\begin{equation}
\label{eq:cmcm+1inati}
\Delta_{c^{m-1} \omega_k, c^{m}\omega_k}(x)
= a^{c^{m-1} \omega_k} \prod_{j=1}^n t_j^{[c^{m-1} \beta_k : \alpha_j]}.
\end{equation}
Our proof of \eqref{eq:cmcm+1inati} uses a little representation theory.
Consider the ring of regular functions $\CC[G]$ as a $G \times G$-representation
under the action $(g_1, g_2)f (x) = f(g_1^T x g_2)$ (see \eqref{eq:T} for the definition of $g_1^T$).
We denote by $f \mapsto (v_1, v_2) f$ the corresponding action of $U(\lg) \times U(\lg)$,
where $U(\lg)$ is the universal enveloping algebra of $\lg$.
For every $f \in \CC[G]$ and $x \in G^{c, c^{-1}}$ as above,~$f(x)$ is a polynomial
in $u_i$ and $t_i$, and the coefficient of
each monomial $u_1^{h_1} \cdots u_n^{h_n}t_1^{d_1} \cdots t_n^{d_n}$
is equal to $((e_1^{(h_1)} \cdots e_n^{(h_n)}, e_{n}^{(d_n)} \cdots e_{1}^{(d_1)})f) (a)$,
where $e_i^{(d)}$ stands for the divided power $e_i^d/d!$ (cf. \cite[Lemma~3.7.5]{bfz}).
Furthermore, if $f$ is a weight vector of weight $(\gamma, \gamma')$ then $f(a)$ can be nonzero
only if $\gamma = \gamma'$.

Now recall that, in view of \eqref{eq:minor-weight}, each
generalized minor $\Delta_{\gamma, \delta}$ with $\gamma, \delta \in W \omega_k$
is a weight vector of weight $(\gamma, \delta)$ in the $G \times G$-irreducible
representation with highest weight $(\omega_k, \omega_k)$.
Recall also the following standard facts from the representation
theory of~$G$: the weight polytope of the fundamental $G$-representation
$V_{\omega_k}$ with highest weight $\omega_k$ has $W \omega_k$ as the set of vertices,
and for each weight $\delta \in W \omega_k$, the corresponding weight subspace
$V_{\omega_k}(\delta)$ is one-dimensional.
Together with basic facts about the representations of $SL_2$,
this implies the following: if $\delta \in W \omega_k$, and $s_i \delta \geq \delta$ for some~$i$
then $s_i \delta - \delta = d_i \alpha_i$, where
$e_i^{(d)}(V_{\omega_k}(\delta)) = \{0\}$ for $d > d_i$;
furthermore, $e_i^{(d_i)}$ establishes the same isomorphism
$V_{\omega_k}(\delta) \to V_{\omega_k}(s_i \delta)$ as the group element ${\overline {s_i}}^{\ -1} \in G$.

Applying all this to $\gamma = c^{m-1} \omega_k$ and $\delta = c^{m} \omega_k$
with $1 \leq m \leq h(k;c)$, we conclude that the only tuple of nonnegative
integers $(h_1, \dots, h_n; d_1, \dots, d_n)$ such that
$(e_1^{(h_1)} \cdots e_n^{(h_n)}, e_{n}^{(d_n)} \cdots e_{1}^{(d_1)})
(\Delta_{\gamma, \delta}) \neq  0$ is the one with all $h_i$ equal
to~$0$, and $\sum_i d_i \alpha_i = \gamma - \delta$.
Furthermore, for this tuple we have
$$((e_1^{(h_1)} \cdots e_n^{(h_n)}, e_{n}^{(d_n)} \cdots e_{1}^{(d_1)})
\Delta_{\gamma, \delta})(a) = \Delta_{\gamma, \gamma}(a) = a^\gamma,$$
implying \eqref{eq:cmcm+1inati}.

We now show that
\begin{equation}
\label{eq:cmcm-2}
\Delta_{c^{m-1} \omega_k, c^{m}\omega_k}(x)
= \prod_{i=1}^n \Delta_{\omega_i, c \omega_i}(x)^{[c^{m-1} \beta_k : \beta_i]}
\end{equation}
for any $x \in G^{c, c^{-1}}$.
It is enough to show that \eqref{eq:cmcm-2} holds for~$x$ of the
form \eqref{eq:Gc-c-1-factor}.
Then both sides of \eqref{eq:cmcm-2} are given by
\eqref{eq:cmcm+1inati}.
Comparing the contributions of $a \in H$ and of each $t_j$ on both
sides, we only have to prove the following two identities:
\begin{equation}
\label{eq:a-contribution}
c^{m-1} \omega_k = \sum_{i=1}^n [c^{m-1} \beta_k : \beta_i] \omega_i;
\end{equation}
\begin{equation}
\label{eq:tj-contribution}
[c^{m-1} \beta_k : \alpha_j] = \sum_{i=1}^n [c^{m-1} \beta_k : \beta_i] [\beta_i : \alpha_j].
\end{equation}

The identity \eqref{eq:tj-contribution} is immediate, so let us
prove \eqref{eq:a-contribution}.
Abbreviating $[c^{m-1} \beta_k : \beta_i] = b_i$ and recalling
\eqref{eq:cm-beta-k}, we can rewrite the equality
$$c^{m-1} \beta_k  =  \sum_{i=1}^n b_i \beta_i$$
in the form
$$c^{m-1} \omega_k  -  \sum_{i=1}^n b_i \omega_i =
c(c^{m-1} \omega_k  -  \sum_{i=1}^n b_i \omega_i).$$
To finish the proof of \eqref{eq:a-contribution}, it remains to use the well-known property
that no Coxeter element has~$1$ as an eigenvalue (see, e.g.,
\cite[Lemma~8.1]{kostant-1}).

In view of \eqref{eq:cmcm-2} and
\eqref{eq:product-coeffs-simplified}, to finish the proof of
Lemma~\ref{lem:coeffs-nonspecial} (hence that of the relation
\eqref{eq:exchange-primitive-nonspecial}), it remains
to show that the restriction to $L^{c,c^{-1}}$ of each
$\Delta_{c^{m} \omega_k, c^{m-1} \omega_k}$ is equal to~$1$.
In view of \eqref{eq:minor-transpose}, the identity \eqref{eq:cmcm-2} implies that
$$\Delta_{c^{m} \omega_k, c^{m-1} \omega_k}(x)
= \prod_{i=1}^n \Delta_{c \omega_i,  \omega_i}(x)^{[c^{m-1} \beta_k : \beta_i]}$$
for any $x \in G^{c, c^{-1}}$.
Now the fact that the right-hand side of this equality is equal to~$1$
for $x \in L^{c, c^{-1}}$, follows from Proposition~\ref{pr:Luv-equations}, and we are done.
\end{proof}

Let us now turn to the proof of \eqref{eq:exchange-primitive-special}.
It is enough to show that both sides of this identity take the
same values at a general element $x \in L^{c,c^{-1}}$ expressed in
the form \eqref{eq:Lc-c-1-factor}.

To compute $\Delta_{\omega_i, \omega_i}(x) = ([x]_0)^{\omega_i}$
(see \eqref{eq:Delta-general}), express each factor $x_{-j}(u_j)$ in $x$ as
$x_{\bar j}(u_j) u_j^{- \alpha_j^\vee}$, and move the
factors $u_j^{- \alpha_j^\vee}$ all the way to the ``center'', using
the commutation relations in \cite[(2.5)]{fz-double}:
\begin{equation}
\label{eq:H-conjugation}
a x_j (t) = x_j (a^{\alpha_j} t) a \ , \quad
a x_{\bar j} (u) = x_{\bar j} (a^{- \alpha_j} u) a \quad (a \in H) \ .
\end{equation}
We obtain $$x = x_{\bar 1}(r_1) \cdots x_{\bar n}(r_n)\ (\prod_{j
= 1}^n u_j^{- \alpha_j^\vee})\ x_n(t_n) \cdots x_1(t_1),$$
where the precise expression for~$r_j$ is immaterial for our current purposes.
It follows that $[x]_0 = \prod_{j = 1}^n  u_j^{- \alpha_j^\vee}$, hence
\begin{equation}
\label{eq:initial-principal-minors}
\Delta_{\omega_i, \omega_i}(x) = \prod_{j = 1}^n  u_j^{\omega_i (- \alpha_j^\vee)} = u_i^{-1}.
\end{equation}

In view of (\ref{eq:cmcm+1inati}), we also have
$$\Delta_{\omega_i, c\omega_i}(x) = u_i^{-1} \prod_{j=1}^n t_j^{[\beta_i : \alpha_j]}.$$
Remembering the definition \eqref{eq:y-j-c} and using (\ref{eq:beta-telescoping}), we get
\begin{equation}
\label{eq:initial-yj}
y_{k;c}(x) =  \Delta_{\omega_k, c \omega_k}(x) \prod_{i < k}
\Delta_{\omega_i, c \omega_i}(x)^{a_{i,k}}
= t_k u_k^{-1} \prod_{j < k} u_j^{-a_{j,k}}.
\end{equation}

Substituting the expressions in \eqref{eq:initial-principal-minors} and \eqref{eq:initial-yj}
into \eqref{eq:exchange-primitive-special}, we can simplify it as follows:
\begin{equation}
\label{eq:exchange-primitive-special-simplified}
\Delta_{-\omega_k, -\omega_k}(x)  =
t_k \prod_{j > k} \Delta_{-\omega_j, -\omega_j}(x)^{-a_{j,k}} + u_k.
\end{equation}

To prove \eqref{eq:exchange-primitive-special-simplified}, note
that, according to (~\ref{eq:minors-tau}), we have
$$\Delta_{-\omega_k, -\omega_k}(x) = \Delta_{\omega_k, \omega_k}(x^\iota).$$
Here $x \mapsto x^\iota$ is an involutive antiautomorphism of the group $G$
defined in \eqref{eq:iota}.
As an easy consequence of this definition, we have
$$x_{-j}(u)^\iota = x_{-j}(u^{-1}),$$
hence
$$x^\iota = x_1(t_1) \cdots x_n(t_n)  x_{-n}(u_n^{-1})  \cdots x_{-1}(u_1^{-1}).$$

To compute $\Delta_{\omega_k, \omega_k}(x^\iota)$ we use the
commutation relations in \cite[Proposition~7.2]{bz01}, which we
rewrite as follows:
$$x_j(t) x_{-i}(u^{-1}) = x_{-i}(u^{-1})x_j(t u^{-a_{i,j}}) \quad
(i \neq j),$$
$$x_j(t) x_{-j}(u^{-1}) = x_{-j}((t+u)^{-1})x_j(w),$$
where the precise expression for~$w$ is immaterial for our purposes.
Iterating these relations, we can rewrite $x^\iota$ in the form
$$x^\iota = x_{-n}(v_n^{-1})  \cdots x_{-1}(v_1^{-1})
x_1(w_1) \cdots x_n(w_n),$$
where the parameters $v_1, \dots, v_n$ satisfy the recurrence relations
$$v_k  = t_k \prod_{j > k} v_j^{-a_{j,k}} + u_k.$$
It remains to observe that $\Delta_{\omega_k, \omega_k}(x^\iota) = v_k$
(cf.~\eqref{eq:initial-principal-minors}), so the last
formula implies \eqref{eq:exchange-primitive-special-simplified}.
This completes the proof of \eqref{eq:exchange-primitive-nonspecial}.

\section{Proof of Theorem~\ref{th:ca-general-1}}
\label{sec:Cluster-algebras}


We start by recalling basic definitions and facts on cluster algebras following
mostly \cite{ca4}.
The definition of a cluster algebra~$\Acal$ starts with introducing
its ground ring.
Let~$\PP$ be a \emph{semifield}, i.e.,
an abelian multiplicative group endowed with a binary operation of
\emph{(auxiliary) addition}~$\oplus$ which is commutative, associative, and
distributive with respect to the multiplication in~$\PP$.
The multiplicative group of~$\PP$ is torsion-free
\cite[Section~5]{ca1}, hence its group ring~$\ZZ\PP$---which will be
used as a \emph{ground ring} for~$\Acal$---is a domain.
Let $\QQ(\PP)$ denote the field of fractions of $\ZZ\PP$.

Every finite family $\{u_j: j \in J\}$ gives rise to a
\emph{tropical semifield} $\Trop (u_j: j \in J)$ defined as follows.
As a multiplicative group, $\Trop (u_j: j \in J)$ is an abelian group
freely generated by the elements $u_j \, (j \in J)$; and
the addition~$\oplus$ in $\Trop (u_j: j \in J)$ is given by
\begin{equation}
\label{eq:tropical-addition}
\prod_j u_j^{a_j} \oplus \prod_j u_j^{b_j} =
\prod_j u_j^{\min (a_j, b_j)} \,.
\end{equation}
Thus, the group ring of $\Trop (u_j: j \in J)$ is the ring of Laurent
polynomials in the variables~$u_j\,$.
A cluster algebra associated to a tropical semifield is said to be
of \emph{geometric type}.

As an \emph{ambient field} for a cluster algebra~$\Acal$ of rank~$n$, we take a field $\Fcal$
isomorphic to the field of rational functions in $n$ independent
variables, with coefficients in~$\QQ (\PP)$.
Note that the definition of $\Fcal$ ignores the auxiliary addition
in~$\PP$.

\begin{definition}
[\emph{Seeds}]
\label{def:seed}
A (labeled) \emph{seed} in~$\Fcal$
is a triple $(\xx, \yy, B)$, where
\begin{itemize}
\item
$\xx = \{x_i: i \in I\}$ is an $n$-tuple of elements of~$\Fcal$
forming a \emph{free generating set} for~$\Fcal$,
that is, the elements $x_i$ for $i \in I$ are algebraically independent over~$\QQ (\PP)$, and $\Fcal =
\QQ (\PP)(\xx)$.
\item
$\yy = \{y_i: i \in I\}$ is an $n$-tuple
of elements of $\PP$.
\item
$B = (b_{i,j})_{i, j \in I}$ is an $n\!\times\! n$ integer matrix
which is \emph{skew-symmetrizable},
that is, $d_i b_{i,j} = - d_j b_{j,i}$ for some positive integers~$d_i \,\, (i \in I)$.
\end{itemize}
We refer to~$\xx$ as the \emph{cluster} of a seed $(\xx, \yy, B)$,
to the tuple~$\yy$ as the \emph{coefficient tuple}, and to the
matrix~$B$ as the \emph{exchange matrix}.
\end{definition}

We will use the notation
$[b]_+ = \max(b,0)$.

\begin{definition}
[\emph{Seed mutations}]
\label{def:seed-mutation}
Let $(\xx, \yy, B)$ be a seed in $\Fcal$,
and $k \in I$.
The \emph{seed mutation} $\mu_k$ in direction~$k$ transforms
$(\xx, \yy, B)$ into the seed
$\mu_k(\xx, \yy, B)=(\xx', \yy', B')$ defined as follows:
\begin{itemize}
\item
The entries of $B'=(b'_{i,j})$ are given by
\begin{equation}
\label{eq:matrix-mutation}
b'_{i,j} =
\begin{cases}
-b_{i,j} & \text{if $i=k$ or $j=k$;} \\[.05in]
b_{i,j} + [b_{ik}]_+ [b_{kj}]_+ - [-b_{i,k}]_+ [-b_{k,j}]_+
 & \text{otherwise.}
\end{cases}
\end{equation}
\item
The coefficient tuple $\yy'$ is given by
\begin{equation}
\label{eq:y-mutation}
y'_j =
\begin{cases}
y_k^{-1} & \text{if $j = k$};\\[.05in]
y_j y_k^{[b_{k,j}]_+}
(y_k \oplus 1)^{- b_{k,j}} & \text{if $j \neq k$}.
\end{cases}
\end{equation}
\item
The cluster $\xx'$ is given by
$x_j'=x_j$ for $j\neq k$,
whereas $x'_k \in \Fcal$ is determined
by the \emph{exchange relation}
\begin{equation}
\label{eq:exchange-rel-xx}
x'_k = \frac
{y_k \ \prod x_i^{[b_{i,k}]_+}
+ \ \prod x_i^{[-b_{i,k}]_+}}{(y_k \oplus 1) x_k} \, .
\end{equation}
\end{itemize}
\end{definition}

It is easy to see that $B'$ is
skew-symmetrizable (with the same choice of the factors $d_i$),
implying that $(\xx', \yy', B')$ is indeed a seed.
Furthermore, the seed mutation $\mu_k$ is involutive,
that is, it transforms $(\xx', \yy', B')$ back into $(\xx, \yy, B)$.
This makes the following equivalence relation on seeds well defined:
we say that two seeds $(\xx, \yy, B)$ and $(\xx', \yy', B')$ are \emph{mutation equivalent}, and write
$(\xx, \yy, B) \sim (\xx', \yy', B')$, if $(\xx', \yy', B')$ can be obtained from $(\xx, \yy, B)$
by a sequence of seed mutations.
For a mutation equivalence class of seeds~$\mathcal{S}$,
we denote by $\mathcal{X} = \mathcal{X}(\mathcal{S})$
the union of clusters of all seeds in~$\mathcal{S}$.
We refer to the elements in $\mathcal{X}$ as \emph{cluster variables}.

\begin{definition}
[\emph{Cluster algebra}]
\label{def:cluster-algebra}
The \emph{cluster algebra} $\Acal$ associated with a mutation
equivalence class of seeds~$\mathcal{S}$
is the $\ZZ \PP$-subalgebra of the ambient field $\Fcal$
generated by all cluster variables: $\Acal = \ZZ \PP[\Xcal]$.
We denote $\Acal = \Acal(\mathcal{S}) = \Acal(\xx, \yy, B) = \Acal(\yy, B)$, where
$(\xx,\yy,B)$ is any seed in~$\mathcal{S}$.
\end{definition}

Following \cite{ca4}, we now introduce an important special system of
coefficients called \emph{principal}.

\begin{definition}
[\emph{Principal coefficients}]
\label{def:principal-coeffs}
We say that a seed $(\xx^\circ, \yy^\circ, B)$ has
\emph{principal coefficients} if the coefficient semifield
$\PP^\circ$ is the tropical semifield $\PP^\circ = \Trop(y_i^\circ: i \in I)$
generated by the coefficient tuple $\yy^\circ$.
\end{definition}

Let $\Acal^\circ = \Acal(\xx^\circ, \yy^\circ, B)$ be the cluster algebra with principal
coefficients at an initial seed $(\xx^\circ, \yy^\circ, B)$.
For convenience, let us use $I = \{1, \dots, n\}$ as the set of
indices for the initial cluster variables $x_1^\circ, \dots, x_n^\circ$ and
the coefficient tuple $y_1^\circ, \dots, y_n^\circ$.
For $j = 1, \dots, n$, we set
\begin{equation}
\label{eq:y-hat}
\widehat y_j^\circ = y_j^\circ \prod_{i=1}^n (x_i^\circ)^{b_{i,j}} \in \Fcal.
\end{equation}
As a special case of \cite[Corollary~6.3]{ca4}, every cluster
variable~$z \in \Acal^\circ$ can be (uniquely) written in the form
\begin{equation}
\label{eq:cluster-vars-principal}
z = (x_1^\circ)^{g_1} \cdots (x_n^\circ)^{g_n} F_{z;\xx^\circ}(\widehat y_1^\circ, \dots, \widehat y_n^\circ)
\end{equation}
for some integer vector
$\gg_{z;\xx^\circ} = (g_1, \dots, g_n) \in \ZZ^n$,
and some integer polynomial
$F_{z;\xx^\circ}(u_1, \dots, u_n) \in \ZZ[u_1, \dots, u_n]$ not
divisible by any variable~$u_i$.
The vector $\gg_{z;\xx^\circ}$ (resp. the polynomial $F_{z;\xx^\circ}$) is
called the \emph{$\gg$-vector} (resp the \emph{$F$-polynomial})
of~$z$ with respect to~$\xx^\circ$.

The $F$-polynomials are conjectured (and in many cases proved)
to have positive coefficients.
However, even with this conjecture not yet proved in complete
generality, one can still evaluate every $F_{z;\xx^\circ}$ at an
$n$-tuple of elements of an arbitrary semifield~$\PP$, since, as shown
in \cite{ca4}, $F_{z;\xx^\circ}(u_1, \dots, u_n)$ can be expressed
as a subtraction-free rational expression in $u_1, \dots, u_n$.
We denote such an evaluation as $F_{z;\xx^\circ}|_\PP$.
Using this notation, one can use $\gg$-vectors and
$F$-polynomials to compute all the cluster variables in an
\emph{arbitrary} cluster algebra
$\Acal = \Acal(\xx, \yy, B)$ with
the coefficients in an arbitrary semifield~$\PP$.
Namely, consider again the cluster algebra $\Acal^\circ = \Acal(\xx^\circ,
\yy^\circ, B)$ with principal coefficients at the initial seed
with the same exchange matrix~$B$.
Then, as shown in \cite[Corollary~6.3]{ca4}, every cluster
variable~$z^\circ \in \Acal^\circ$ gives rise to a cluster variable
$z \in \Acal$ given by
\begin{equation}
\label{eq:cluster-vars-general}
z = x_1^{g_1} \cdots x_n^{g_n}
\frac{F_{z^\circ;\xx^\circ}(\widehat {y}_1, \dots, \widehat {y}_n)}
{F_{z^\circ;\xx^\circ}|_{\PP}({y}_1, \dots, {y}_n)},
\end{equation}
and all the cluster variables in $\Acal$ are of this
form; here the elements $\widehat {y}_j$ of the ambient
field for $\Acal$ have the same meaning as in
\eqref{eq:y-hat}.
Comparing \eqref{eq:cluster-vars-general} and
\eqref{eq:cluster-vars-principal}, we obtain the following
assertion.

\begin{proposition}
\label{pr:principal-embedding}
In the above notation, suppose that the elements ${y}_1, \dots, {y}_n \in \PP$ are
multiplicatively independent, i.e., the
correspondence $y_j^\circ \mapsto y_j$ identifies $\PP^\circ$ with
a multiplicative subgroup of $\PP$.
Then the cluster algebra $\Acal^\circ$ can be identified with a $\ZZ \PP^\circ$-subalgebra of
$\Acal$ via the correspondence sending each cluster
variable~$z^\circ \in \Acal^\circ$ to
$F_{z^\circ;\xx^\circ}|_{\PP}({y}_1, \dots, {y}_n) z \in \Acal$.
With this identification, $\Acal$ is obtained from
$\Acal^\circ$ by the extension of scalars from $\ZZ \PP^\circ$ to $\ZZ \PP$.
\end{proposition}

\begin{definition}
\label{def:finite-type}
A cluster algebra $\Acal = \Acal(\mathcal{S})$ is of
\emph{finite type} if the mutation equivalence class~$\mathcal{S}$
consists of finitely many seeds.
\end{definition}

A classification of cluster algebras of finite type was given in \cite{ca2}.
We present it in the form convenient for our current purposes.

\begin{theorem}[\cite{ca2}]
\label{th:fin-type-class}
A cluster algebra $\Acal$ is of finite type
if and only if the exchange matrix at some
seed of~$\Acal$ is of the form~$B(c)$ (see \eqref{eq:Bc})
for some Coxeter element in the Weyl group associated to a Cartan
matrix~$A$ of finite type.
Furthermore, the type of~$A$ in the Cartan-Killing nomenclature
is uniquely determined by the mutation equivalence class of seeds
in~$\Acal$, and, for a given~$A$, all the matrices $B(c)$
associated to different Coxeter elements are mutation equivalent
to each other.
\end{theorem}

By Theorem~\ref{th:fin-type-class}, the property of a cluster
algebra~$\Acal$ to be of finite type does not depend on the choice of a coefficient system.
As shown in \cite{ca2}, the same is true for the structure of the
set~$\Xcal$ of cluster variables and its division into clusters.
To be more precise, fix a Cartan matrix~$A$ of the same type as $\Acal$, and choose an initial seed
$(\xx, \yy, B)$ with the exchange matrix $B = B(t)$, where $t = t_+ t_-$
is the bipartite Coxeter element (see \eqref{eq:coxeter-bipartite}).
As before, let $\Phi$ be the root system associated with~$A$.
Then the cluster variables in~$\Acal$ are in a natural bijection $\alpha \mapsto x[\alpha]$
with the set $\Phi_{\geq -1}$ of ``almost positive" roots (the union of the set of positive roots
and the set of negative simple roots $-\alpha_i$).
(Recall that, for each $\alpha \in \Phi_{\geq -1}$,
the integer vector in~$\ZZ^I$ with the components $d_i =
[\alpha: \alpha_i]$ is the \emph{denominator vector} in the
Laurent expansion of the cluster variable $x[\alpha]$ with respect
to the initial cluster.)
In particular, each initial cluster variable takes the form $x_i = x[-\alpha_i]$.

To describe the clusters in $\Acal$, we recall the two
involutive permutations $\tau_+, \tau_-$ of $\Phi_{\geq -1}$ given by
\begin{equation}
\label{eq:tau-action}
\tau_\varepsilon (\alpha) =
\begin{cases}
\alpha & \text{if $\alpha = -\alpha_j$ with $\varepsilon(j) = -\varepsilon$;} \\[.05in]
t_\varepsilon (\alpha)
 & \text{otherwise.}
\end{cases}
\end{equation}
As shown in \cite[Section~3.1]{yga}, for any
$\alpha, \beta \in \Phi_{\geq -1}$, there is
a well-defined nonnegative integer $(\alpha \| \beta)$
(\emph{compatibility degree}) uniquely
characterized by the following two properties:
\begin{eqnarray}
\label{eq:compatibility-1}
&(- \alpha_i \| \beta) = [[\beta: \alpha_i]]_+ , &
\\
\label{eq:compatibility-2}
&(\tau_\varepsilon \alpha \| \tau_\varepsilon \beta) = (\alpha \|
\beta) , &
\end{eqnarray}
for any $\alpha, \beta \in \Phi_{\geq -1}$, any $i\in I$,
and any sign~$\varepsilon$.
We say that $\alpha$ and $\beta$ are \emph{compatible} if
$(\alpha \| \beta) = 0$ (this is in fact a symmetric relation).
With this terminology in place, the following was shown in \cite{ca2}.

\begin{proposition}
\label{pr:root-clusters}
Under the above parameterization of cluster variables by almost positive roots,
the clusters in~$\Acal$ are exactly the families $\xx(C) = \{x[\alpha]:
\alpha \in C\}$, where $C$ runs over all maximal by inclusion
subsets of $\Phi_{\geq -1}$ consisting of mutually compatible roots.
\end{proposition}

Note that every cluster $\xx(C)$ is included in a unique seed
$(\xx(C), \yy(C), B(C))$ in $\Acal$.
To emphasize the dependence on~$\Acal$, we use the notation
$(\xx(C), \yy(C), B(C)) = (\xx(C)^{(\Acal)}, \yy(C)^{(\Acal)}, B(C))$;
note that the exchange matrix $B(C)$ is independent of $\Acal$,
and it is explicitly described in \cite[Definition~4.2]{ca2}.

Now we introduce a tool to relate to each other two
cluster algebras of the same finite type but with different
coefficient systems.
This is a version of the \emph{coefficient specialization}
discussed in \cite[Section~12]{ca4}.
We will deal only with cluster algebras of geometric type.
Let us start with a simple general proposition.

\begin{proposition}
\label{pr:trop-semifield-homs}
Let $\PP = \Trop(u_{j}: j  \in J)$ and
$\widetilde \PP = \Trop(u_{\widetilde j}: \widetilde j  \in \widetilde J)$ be
tropical semifields, and $\varphi: \widetilde \PP \to \PP$ a
homomorphism of multiplicative groups.
Then $\varphi$ is a homomorphism of semifields if and only if
each $\varphi(u_{\widetilde j}) \in \PP$ is a monomial with nonnegative exponents
in the generators $u_j$ of $\PP$, and no two such monomials have a
common factor $u_j$ for some $j \in J$.
\end{proposition}

\begin{proof}
This follows at once from the equalities $u_{\widetilde j} \oplus 1
= 1$ and $u_{\widetilde j} \oplus u_{\widetilde {j'}} = 1$ for
$\widetilde j \neq \widetilde {j'}$.
\end{proof}

Now let $\widetilde \Acal$ be a cluster algebra of finite
type as above, with the coefficient semifield
$\widetilde \PP = \Trop(u_{\widetilde j}: \widetilde j  \in \widetilde J)$.
In the above notation, the seeds of $\widetilde \Acal$ are of the form
$(\xx(C)^{(\widetilde \Acal)}, \yy(C)^{(\widetilde \Acal)}, B(C))$.
Let $\varphi: \widetilde \PP \to \PP$ be a homomorphism of tropical semifields
as in Proposition~\ref{pr:trop-semifield-homs}.
Using~$\varphi$, we can form a cluster algebra $\Acal$ of the same
type as $\widetilde \Acal$, but with coefficients in $\PP$: its
seeds are of the form $(\xx(C)^{(\Acal)}, \yy(C)^{(\Acal)},
B(C))$, where the coefficient tuple $\yy(C)^{(\Acal)}$ at some (equivalently, any) seed
is taken as $\varphi(\yy(C)^{(\widetilde \Acal)})$.
The following relationship between $\Acal$ and $\widetilde \Acal$ is
immediate from the definitions.

\begin{proposition}
\label{pr:coef-specialization}
The semifield homomorphism~$\varphi$ induces an isomorphism of $\ZZ \PP$-algebras
$$\varphi_\star: \widetilde \Acal \otimes_{\ZZ \widetilde \PP} \ZZ \PP \to \Acal$$
given by $\varphi_\star (x[\alpha]^{(\widetilde \Acal)} \otimes 1) =
x[\alpha]^{(\Acal)}$ for $\alpha \in \Phi_{\geq -1}$.
\end{proposition}

Turning to Theorem~\ref{th:ca-general-1}, we will show that
Proposition~\ref{pr:coef-specialization} applies to
$\widetilde \Acal_\CC = \CC[G^{c,c^{-1}}]$ and
$\Acal_\CC = \CC[L^{c,c^{-1}}]$, where $\widetilde \Acal_\CC$ and
$\Acal_\CC$ are obtained from $\widetilde \Acal$ and
$\Acal$ by extension of scalars from $\ZZ$ to $\CC$.
First of all, the fact that the coordinate ring of $G^{c,c^{-1}}$
carries the structure of a (complexified) geometric type cluster
algebra $\widetilde \Acal_\CC$ of the same finite type as $G$, was established in \cite[Example~2.24]{ca3}.
We will use a slightly modified description of its initial seed by
choosing a different reduced word for $(c, c^{-1})$.
Namely, we have:
\begin{itemize}
\item
the ambient field
for $\widetilde \Acal_\CC$ is
just the field of rational functions on $G^{c,c^{-1}}$;
\item the coefficient semifield $\widetilde \PP$ is the tropical semifield
with the $2n$ generators $\Delta_{c\omega_i, \omega_i}|_{G^{c,c^{-1}}}$
and $\Delta_{\omega_i, c \omega_i}|_{G^{c,c^{-1}}}$ for $i \in I$;
\item
the initial cluster consists of the functions
$\Delta_{\omega_i, \omega_i}|_{G^{c,c^{-1}}}$;
\item the initial coefficient tuple consists of the restrictions
to $G^{c,c^{-1}}$ of the functions
$\Delta_{\omega_j, c \omega_j} \Delta_{c \omega_j, \omega_j} \prod_{i \prec_c j}
(\Delta_{\omega_i, c \omega_i}\Delta_{c \omega_i, \omega_i})^{a_{i,j}}$;
\item the initial exchange matrix is $B(c)$.
\end{itemize}

Now recall that by Proposition~\ref{pr:Luv-equations},
the reduced double cell $L^{c,c^{-1}}$ is obtained from
$G^{c,c^{-1}}$ by specializing each invertible regular function
$\Delta_{c \omega_i, \omega_i}$ to~$1$.
Algebraically, this means that the coordinate ring
$\CC[L^{c,c^{-1}}]$ can be described as follows.
Let $\PP$ be the tropical semifield
with the $n$ generators $\Delta_{\omega_i, c \omega_i}|_{L^{c,c^{-1}}}$
for $i \in I$, and let $\varphi: \widetilde \PP \to \PP$ be the
restriction homomorphism from $G^{c,c^{-1}}$ to $L^{c,c^{-1}}$.
Then $\varphi$ is as in Proposition~\ref{pr:trop-semifield-homs},
acting on the generators by
$$\varphi(\Delta_{c\omega_i, \omega_i}|_{G^{c,c^{-1}}}) = 1,
\quad \varphi(\Delta_{\omega_i, c \omega_i}|_{G^{c,c^{-1}}}) =
\Delta_{\omega_i, c \omega_i}|_{L^{c,c^{-1}}};$$
and we have
$$\CC[L^{c,c^{-1}}] = \CC[G^{c,c^{-1}}] \otimes_{\ZZ \widetilde \PP} \ZZ \PP.$$
Comparing this with Proposition~\ref{pr:coef-specialization}, we
conclude that $\CC[L^{c,c^{-1}}]$ is the complexified cluster
algebra $\Acal_\CC$ with the coefficient semifield $\PP$ and the initial seed
as in Theorem~\ref{th:ca-general-1}.

This is still not quite what we need to prove
Theorem~\ref{th:ca-general-1}, since the auxiliary addition in
$\PP$ is \emph{different} from the one in
$\PP^\circ = \Trop(y_{j;c}: j  \in I)$.
However, the elements $y_{j;c} \in \PP$ are related with
the generators $\Delta_{\omega_i, c \omega_i}|_{L^{c,c^{-1}}}$ by
a triangular (hence invertible) monomial transformation.
Therefore, as a multiplicative group, $\PP^\circ$ coincides with~$\PP$.
By Proposition~\ref{pr:principal-embedding}, the cluster
algebra $\Acal^\circ$ with principal coefficients at the initial
seed can be identified with $\Acal$ (with the cluster variables in $\Acal^\circ$
obtained by rescaling from those in $\Acal$).
This concludes the proof of Theorem~\ref{th:ca-general-1}.

\section{Proofs of Theorems~\ref{th:ca-general-2} and~\ref{th:ca-general-3}}
\label{sec:univcoeff}

In this section we fix an indecomposable Cartan matrix~$A$ of finite
type, and freely use the root system formalism developed in
Section~\ref{sec:cce}, and the cluster algebra formalism developed in
Section~\ref{sec:Cluster-algebras} (as before, lifting the restriction that $A$
is indecomposable presents no problem).

Let us recap a little.
For a Coxeter element~$c$ in the Weyl group of~$A$, let
$\Acal^\circ(c)$ be the cluster algebra with principal
coefficients at an initial seed $(\xx^\circ, \yy^\circ, B(c))$ as
defined in Definition~\ref{def:principal-coeffs}.
We have already shown in Theorem~\ref{th:ca-general-1} that the
complexification of $\Acal^\circ(c)$ can be identified with
the coordinate ring $\CC[L^{c,c^{-1}}]$, and under this
identification, the initial cluster variables $x_i^\circ = x_{\omega_i;c}$ and the initial
coefficients $y_j^\circ = y_{j;c}$  are given by \eqref{eq:x-gamma-c}
and \eqref{eq:y-j-c}, respectively.
We have also proved that the elements $x_{\gamma;c}$
(for $\gamma \in \Pi(c)$) and $y_{j;c}$ satisfy algebraic relations
\eqref{eq:exchange-primitive-special} and
\eqref{eq:exchange-primitive-nonspecial} in $\Acal^\circ(c)$.
Thus, to prove Theorem~\ref{th:ca-general-2} and
to complete the proof of Theorem~\ref{th:ca-general-3}, we can
forget about the geometric realization of $\Acal^\circ(c)$, and it remains to
prove the following two statements:
\begin{align}
\label{eq:cluster-vars-principal-list}
&\text{The elements $x_{\gamma;c} \in \Acal^\circ(c)$
(for $\gamma \in \Pi(c)$) satisfying \eqref{eq:exchange-primitive-special}}\\
\nonumber
&\text{and \eqref{eq:exchange-primitive-nonspecial} are exactly the cluster
variables in  $\Acal^\circ(c)$.}
\end{align}
\begin{align}
\label{eq:primitive-exchange-principal}
&\text{The relations \eqref{eq:exchange-primitive-special} and
\eqref{eq:exchange-primitive-nonspecial} are exactly}\\
\nonumber
&\text{the primitive
exchange relations in  $\Acal^\circ(c)$.}
\end{align}

For \eqref{eq:primitive-exchange-principal} recall from \cite[Section~12]{ca4} that an
exchange relation \eqref{eq:exchange-rel-xx} is called
\emph{primitive} if one of the products of cluster variables in
the right hand side is empty, understood to be equal to~$1$.

We start with some combinatorial preparation.
First we transfer the compatibility degree function
given by \eqref{eq:compatibility-1} and \eqref{eq:compatibility-2}
from $\Phi_{\geq -1}$ to each $\Pi(c)$.
Let $\tau_c$ be the permutation of $\Pi(c)$ given by
\begin{equation}
\label{eq:tau-c}
\tau_c(c^{m-1} \omega_i) = c^{m} \omega_i \quad (1 \leq m \leq h(i;c)),
\quad \tau_c(- \omega_i) = \omega_i.
\end{equation}

\begin{proposition}
\label{pr:comp-degree-Pic}
There is a unique assignment of a nonnegative integer $(\gamma \| \delta)_c$
to any pair of weights $\gamma, \delta \in \Pi(c)$, satisfying the following properties:
\begin{enumerate}
\item $(\omega_i \| \omega_j)_c = 0 \quad (i, j \in I)$,
\item $(\omega_i \| \gamma)_c = [c^{-1} \gamma -  \gamma  : \alpha_i]
\quad (\gamma \in \Pi(c) - \{\omega_j: j \in I\})$,
\item $(\tau_c \gamma \| \tau_c \delta)_c = (\gamma \| \delta)_c \quad
(\gamma, \delta \in \Pi(c))$.
\end{enumerate}
\end{proposition}

We call $(\gamma \| \delta)_c$ the
\emph{$c$-compatibility degree} of $\gamma$ and $\delta$.

\begin{proof}
The uniqueness of the $c$-compatibility degree with desired properties is clear.
Note also that, in view of \eqref{eq:omega-c-omega},
condition (2) in Proposition~\ref{pr:comp-degree-Pic} can be rewritten as
\begin{equation}
\label{eq:c-comp-degree-2}
(\omega_i \| c^m \omega_j)_c = [c^{m-1} \beta_j  : \alpha_i]
\quad (i, j \in I; \, 1 \leq m \leq h(j;c)),
\end{equation}
where the root $\beta_j$ is given by \eqref{eq:beta-i}.

In view of \eqref{eq:cyclical-moves-transitive}, to prove the existence it suffices to
do it for the bipartite Coxeter element $t = t_+
t_-$, and then to show that the existence of the $c$-compatibility degree
with desired properties implies that of $\tilde c$-compatibility
degree, where $\tilde c$ is obtained from $c$
as in \eqref{eq:cyclical-move}.

We start by dealing with $\Pi(t)$.
Consider the permutation  $\tau = \tau_+ \tau_-$ of
$\Phi_{\geq -1}$, where $\tau_+$ and $\tau_-$ are given by \eqref{eq:tau-action}.

\begin{lemma}
\label{lem:bijection-almost-positive-roots-Pit}
There is a bijection $\psi: \Pi(t) \to \Phi_{\geq -1}$ uniquely
determined by the properties that $\psi(\omega_i) = - \alpha_i$
for all $i \in I$, and $\psi \circ \tau_t = \tau \circ \psi$.
Explicitly, $\psi$ is given by
\begin{equation}
\label{eq:psi}
\psi(\gamma) =
\begin{cases}
-\alpha_i & \text{if $\gamma = \omega_i$ for some $i \in I$;} \\[.1in]
t^{-1} \gamma - \gamma & \text{otherwise.}
\end{cases}
\end{equation}
\end{lemma}

\begin{proof}
The uniqueness of~$\psi$ with the desired properties is clear.
The fact that \eqref{eq:psi} defines a bijection between $\Pi(t)$ and
$\Phi_{\geq -1}$ follows from Lemmas~\ref{lem:two-identities} and
\ref{lem:c-reps-in-PHI}.
It remains to prove that the map $\psi$ given by \eqref{eq:psi}
satisfies the property that $\psi(\tau_t \gamma) = \tau
\psi(\gamma)$ for $\gamma \in \Pi(t)$.
This is clear if both $\psi(\gamma)$ and $\psi(\tau_t \gamma)$
fall into the second case in \eqref{eq:psi}, i.e., if $\gamma$ is
not of the form $\pm \omega_i$.
By the definitions, we also have
\begin{equation*}
\psi(\tau_t \omega_i) = \tau
\psi(\omega_i) =
\begin{cases}
\alpha_i & \text{if $\varepsilon(i) = +1$;} \\[.1in]
t_+ \alpha_i  & \text{if $\varepsilon(i) = -1$,}
\end{cases}
\end{equation*}
and $\psi(\tau_t (-\omega_i)) = \tau
\psi(-\omega_i) = -\alpha_i$ for all $i \in I$, finishing the
proof.
\end{proof}

By Lemma~\ref{lem:bijection-almost-positive-roots-Pit}, one
can define the $t$-compatibility degree by setting
\begin{equation}
\label{eq:t-comp}
(\gamma \| \delta)_t = (\psi(\gamma) \| \psi(\delta))
\quad (\gamma, \delta \in \Pi(t)),
\end{equation}
and  it satisfies the desired properties.

Now suppose that $\tilde c$ is obtained from $c$ via \eqref{eq:cyclical-move}.
Without loss of generality, we assume that $I = \{1, \dots, n\}$,
$c = s_1 s_2 \cdots s_n$, and $\tilde c = s_2 \cdots s_n s_1$.
It remains to show that the
existence of the $c$-compatibility degree with
the desired properties implies that of the $\tilde c$-compatibility degree.
The following lemma is an easy consequence of \eqref{eq:h(i;tildec)}.

\begin{lemma}
\label{lem:bijection-Pic-Pitildec}
There is a bijection $\psi_{c,\tilde c}: \Pi(\tilde c) \to \Pi(c)$ uniquely
determined by the properties that
\begin{equation}
\label{eq:psi-tilde-c}
\psi_{c,\tilde c}(\omega_i) =
\begin{cases}
\omega_i & \text{if $i \neq 1$;} \\[.1in]
c \omega_i & \text{if $i = 1$,}
\end{cases}
\end{equation}
and $\psi_{c,\tilde c} \circ \tau_{\tilde c} = \tau_c \circ \psi_{c,\tilde c}$.
Explicitly, for $\gamma \in \Pi(\tilde c)$ we have
\begin{equation}
\label{eq:psi-c-tildec-explicit}
\psi_{c,\tilde c}(\gamma) =
\begin{cases}
\omega_1 & \text{if $\gamma = -\omega_1$;} \\[.1in]
s_1 \gamma & \text{otherwise.}
\end{cases}
\end{equation}
\end{lemma}

Now we define the $\tilde c$-compatibility degree by setting
\begin{equation}
\label{eq:tilde-c-thru-c-comp}
(\gamma \| \delta)_{\tilde c} = (\psi_{c,\tilde c}(\gamma) \| \psi_{c,\tilde c}(\delta))_c
\quad (\gamma, \delta \in \Pi(\tilde c)).
\end{equation}
This makes condition (3) in Proposition~\ref{pr:comp-degree-Pic} obvious.
It remains to check conditions (1) and (2).

The equality $(\omega_i \| \omega_j)_{\tilde c} = 0$ is obvious if $i, j \neq 1$.
Let us show that $(\omega_1 \| \omega_i)_{\tilde c} = 0$ and $(\omega_i \| \omega_1)_{\tilde c} = 0$
for $i \neq 1$.
Indeed, we have
\begin{align*}
(\omega_1 \| \omega_i)_{\tilde c} &= (c\omega_1 \| \omega_i)_{c}
= (\omega_1 \| -\omega_i)_{c}\\
&= [\omega_i - c^{-1} \omega_i  :
\alpha_1] = [s_n s_{n-1} \cdots s_{i+1} \alpha_i  :  \alpha_1] = 0,
\end{align*}
and
$$(\omega_i \| \omega_1)_{\tilde c} = (\omega_i \| c \omega_1)_{c}
= [\omega_i - c \omega_1  :  \alpha_i] = [\alpha_1  :  \alpha_i] = 0.$$

It remains to prove that $(\omega_i \| \tilde c^m \omega_j)_{\tilde c}$
for $i, j \in \{1, \dots, n\}$ and $1 \leq m \leq h(j;\tilde c)$
is given by \eqref{eq:c-comp-degree-2}, i.e., we have
\begin{equation}
\label{eq:tilde-c-comp-2}
(\omega_i \| \tilde c^m \omega_j)_{\tilde c} =
[\tilde c^{m-1} s_2 \cdots s_{j-1} \alpha_j  : \alpha_i],
\end{equation}
with the convention that for $j=1$, the index $j-1$ is understood to
be equal to~$n$.

In checking \eqref{eq:tilde-c-comp-2}, we will repeatedly use the following obvious
property: the coefficient $[\alpha  :  \alpha_i]$ does not
change under replacing $\alpha$ with $s_k \alpha$ for $k \neq i$.
Now let us consider four separate cases.

\smallskip

{\bf Case 1.} Let $i, j \neq 1$.
Then we have
\begin{align*}
(\omega_i \| \tilde c^m \omega_j)_{\tilde c} &= (\omega_i \| c^m \omega_j)_{c}
= [c^{m-1} s_1 \cdots s_{j-1} \alpha_j  : \alpha_i]\\
&= [s_1 \tilde c^{m-1} s_2 \cdots s_{j-1} \alpha_j  : \alpha_i]
= [\tilde c^{m-1} s_2 \cdots s_{j-1} \alpha_j  : \alpha_i],
\end{align*}
proving \eqref{eq:tilde-c-comp-2}.

\smallskip

{\bf Case 2.} Let $i=1$ and $j \neq 1$.
If $m \geq 2$ then we have
\begin{align*}
(\omega_1 \| \tilde c^m \omega_j)_{\tilde c} &= (\omega_1 \| c^{m-1} \omega_j)_{c}
= [c^{m-2} s_1 \cdots s_{j-1} \alpha_j  : \alpha_1]\\
&= [s_2 \cdots s_n c^{m-2} s_1 \cdots s_{j-1} \alpha_j  : \alpha_1]
= [\tilde c^{m-1} s_2 \cdots s_{j-1} \alpha_j  : \alpha_1],
\end{align*}
proving \eqref{eq:tilde-c-comp-2}.
And if $m =1$ then
$$(\omega_1 \| \tilde c \omega_j)_{\tilde c} = (\omega_1 \|
\omega_j)_{c}= 0 = [s_2 \cdots s_{j-1} \alpha_j  : \alpha_1],$$
again proving \eqref{eq:tilde-c-comp-2}.

\smallskip

{\bf Case 3.} Let $i \neq 1$ and $j = 1$.
Then we have
\begin{align*}
(\omega_i \| \tilde c^m \omega_1)_{\tilde c} &=
(\omega_i \| c^{m+1} \omega_1)_{c} =
[c^{m} \alpha_1  : \alpha_i]\\ & =
[s_1 \tilde c^{m-1} s_2 \cdots s_{n} \alpha_1  : \alpha_i] =
[\tilde c^{m-1} s_2 \cdots s_{n} \alpha_1  : \alpha_i],
\end{align*}
proving \eqref{eq:tilde-c-comp-2}.

\smallskip

{\bf Case 4.} Let $i = j = 1$.
Then we have
\begin{align*}
(\omega_1 \| \tilde c^m \omega_1)_{\tilde c} & =
(\omega_1 \| c^{m} \omega_1)_{c} = [c^{m-1} \alpha_1  : \alpha_1]\\& =
[s_2 \cdots s_{n} c^{m-1} \alpha_1  : \alpha_1] =
[\tilde c^{m-1} s_2 \cdots s_{n} \alpha_1  : \alpha_1],
\end{align*}
proving \eqref{eq:tilde-c-comp-2} in this case as well.

\smallskip
This completes the proof of Proposition~\ref{pr:comp-degree-Pic}.
\end{proof}

We now show that the $c$-compatibility degree satisfies ``Langlands duality."
Let $\Phi^\vee$ denote the dual root system to $\Phi$; it
corresponds to the transpose Cartan matrix~$A^T$.
The Weyl group of $A^T$ is identified with the Weyl
group~$W$ of~$A$, and there is a canonical $W$-equivariant bijection $\alpha
\mapsto \alpha^\vee$ between $\Phi$ and $\Phi^\vee$.
This bijection restricts to a bijection between $\Phi_{\geq -1}$ and $\Phi^\vee_{\geq -1}$.
As shown in \cite[Proposition~3.3]{yga}, we have
\begin{equation}
\label{eq:root-comp-duality}
(\alpha \| \beta) = (\beta^\vee \| \alpha^\vee) \quad (\alpha, \beta \in \Phi_{\geq -1}).
\end{equation}

Now consider the set $\Pi(c)^\vee$ defined in the same way as
$\Pi(c)$ but associated with the dual root system.
Thus, the elements of $\Pi(c)^\vee$ are coweights of the form
$c^m \omega_i^\vee$, where the $\omega_i^\vee$ are fundamental
coweights, and $0 \leq m \leq h(i;c)$ (it easily follows from
Lemmas~\ref{lem:two-identities} and \ref{lem:c-reps-in-PHI} that
the numbers $h(i;c)$ for the dual root system $\Phi^\vee$ are the same as for $\Phi$).
The correspondence $\omega_i \mapsto \omega_i^\vee$ uniquely extends to
the $\tau_c$-equivariant bijection $\gamma \mapsto \gamma^\vee$
between $\Pi(c)$ and $\Pi(c)^\vee$.
By the definition, for every $\gamma \in \Pi(c) - \{\omega_i  :  i \in I\}$, we have
\begin{equation}
\label{eq:distance-duality}
(\tau_c^{-1} \gamma - \gamma)^\vee = \tau_c^{-1} \gamma^\vee -
\gamma^\vee.
\end{equation}

With some abuse of notation, we denote by $(\gamma^\vee \| \delta^\vee)_c$ the $c$-compatibility degree
on $\Pi(c)^\vee$.

\begin{proposition}
\label{pr:comp-duality-c}
For $\gamma, \delta \in \Pi(c)$, we have
\begin{equation}
\label{eq:comp-duality-c}
(\gamma \| \delta)_c = (\delta^\vee \| \gamma^\vee)_c.
\end{equation}
\end{proposition}

\begin{proof}
The equality \eqref{eq:comp-duality-c} follows at once from
\eqref{eq:root-comp-duality} and an obvious fact that each of the
bijections $\psi$ and $\psi_{c,\tilde c}$ (see
Lemmas~\ref{lem:bijection-almost-positive-roots-Pit} and
\ref{lem:bijection-Pic-Pitildec}) commutes with passing to the
dual roots and weights.
\end{proof}

We illustrate the use of \eqref{eq:comp-duality-c} by the
following lemma to be used later.

\begin{lemma}
\label{lem:comp-c-identity}
For any index $j \in I$, and $\gamma \in \Pi(c) - \{\omega_j, c \omega_j\}$, we have
\begin{equation}
\label{eq:comp-c-omega-omega}
(\gamma \| c \omega_j)_c + \sum_{i \prec_c j} a_{i,j} (\gamma \| c \omega_i)_c =
-((\gamma \| \omega_j)_c + \sum_{j \prec_c i} a_{i,j} (\gamma \| \omega_i)_c).
\end{equation}
Furthermore, for $\gamma = \omega_j$, the right hand side of \eqref{eq:comp-c-omega-omega} is equal
to~$0$, while the left hand side is equal to~$1$;
and for $\gamma = c \omega_j$,
the left hand side of \eqref{eq:comp-c-omega-omega} is equal
to~$0$, while the right hand side is equal to~$-1$.
\end{lemma}

\begin{proof}
First of all, if $\gamma = \omega_k$ for some $k$, then the right hand side of \eqref{eq:comp-c-omega-omega} is equal
to~$0$ by condition (1) in Proposition~\ref{pr:comp-degree-Pic}.
As for the left hand side, by condition (2) it is equal to
$$[\beta_j + \sum_{i \prec_c j} a_{i,j} \beta_i  :  \alpha_k] =
[\alpha_j  :  \alpha_k] = \delta_{jk}$$
(see \eqref{eq:beta-telescoping}).
The case $\gamma = c \omega_k$ is treated similarly: now the left hand side
of \eqref{eq:comp-c-omega-omega} is equal to~$0$, while the right hand side
is equal to $-\delta_{jk}$ (by applying \eqref{eq:beta-telescoping} with
$c$ replaced by $c^{-1}$).

It remains to consider the case when $\gamma$ is not of the form
$\omega_k$ or $c \omega_k$.
Let $\alpha = c^{-1} \gamma - \gamma$.
Using \eqref{eq:comp-duality-c}, we can rewrite the right hand
side of \eqref{eq:comp-c-omega-omega} as
$$- ([\alpha^\vee  :  \alpha_j^\vee] +
\sum_{j \prec_c i} a_{i,j} [\alpha^\vee  :  \alpha_i^\vee]) =
- \langle \alpha^\vee, \omega_j + \sum_{j \prec_c i} a_{i,j}
\omega_i \rangle,$$
where $\langle ?, ? \rangle$ is the standard $W$-invariant pairing between
the coweights and weights.
Similarly, the left hand side of \eqref{eq:comp-c-omega-omega} can
be rewritten as
$$\langle c^{-1} \alpha^\vee, \omega_j + \sum_{i \prec_c j} a_{i,j}
\omega_i \rangle = \langle \alpha^\vee, c(\omega_j + \sum_{i \prec_c j} a_{i,j}
\omega_i) \rangle.$$
It remains to check the identity
$$c(\omega_j + \sum_{i \prec_c j} a_{i,j} \omega_i) =
-(\omega_j + \sum_{j \prec_c i} a_{i,j} \omega_i).$$
Using \eqref{eq:alpha-omega}, \eqref{eq:omega-c-omega} and
\eqref{eq:beta-telescoping}, we obtain
\begin{align*}
-(\omega_j + \sum_{j \prec_c i} a_{i,j} \omega_i) &=
- \alpha_j + \omega_j + \sum_{i \prec_c j} a_{i,j} \omega_i\\
&= (\omega_j - \beta_j) + \sum_{i \prec_c j} a_{i,j} (\omega_i -
\beta_i) = c(\omega_j + \sum_{i \prec_c j} a_{i,j} \omega_i),
\end{align*}
finishing the proof.
\end{proof}

Returning to the proofs of \eqref{eq:cluster-vars-principal-list}
and \eqref{eq:primitive-exchange-principal}, we will deduce these
statements from the following proposition of independent interest.

\begin{proposition}
\label{pr:universal-ca-c}
For every Coxeter element~$c$, there exists a cluster algebra
$\widetilde \Acal(c)$ of geometric type, satisfying the
following properties:
\begin{enumerate}
\item The coefficient semifield of $\widetilde \Acal(c)$ is
$\widetilde \PP(c) = \Trop (p[\gamma]:  \gamma \in \Pi(c))$.
\item The cluster variables in $\widetilde \Acal(c)$ are labeled
by the set $\Pi(c)$, the variable corresponding to $\gamma \in
\Pi(c)$ being denoted $x[\gamma]$.
\item The initial seed of $\widetilde \Acal(c)$ is of the form
$(\xx, \yy, B(c))$, where $\xx = (x[\omega_i]:  i \in I)$,
and $\yy = (y_j:  j \in I)$ with
\begin{equation}
\label{eq:y-initial-c-universal}
y_j = p[\omega_j]  p[c \omega_j]^{-1}
\prod_{\gamma \in \Pi(c) - \{\omega_j, c \omega_j\}}
p[\gamma]^{(\gamma \| c \omega_j)_c + \sum_{i \prec_c j} a_{i,j} (\gamma \| c
\omega_i)_c}.
\end{equation}
\item The primitive exchange relations in $\widetilde \Acal(c)$
are exactly the following (for $k \in I$ and $1 \leq m \leq h(k;c)$):
\begin{equation}
\label{eq:exchange-primitive-special-universal}
x[-\omega_k]\ x[\omega_k] =
p[\omega_k] \prod_{i \prec_c k} x[\omega_i]^{-a_{i,k}}
\prod_{k \prec_c i} x[-\omega_i]^{-a_{i,k}} + \prod_{\gamma \in \Pi(c)}
p[\gamma]^{(\gamma \| \omega_k)_c};
\end{equation}
\begin{align}
\label{eq:exchange-primitive-nonspecial-universal}
x[c^{m-1} \omega_k] \ x[c^{m} \omega_k] &=
p[c^{m} \omega_k]\prod_{i \prec_c k} x[c^{m} \omega_i]^{-a_{i,k}}
\prod_{k \prec_c i} x[c^{m-1} \omega_i]^{-a_{i,k}}\\
\nonumber
& +
\prod_{\gamma \in \Pi(c)}
p[\gamma]^{(\gamma \| c^{m}\omega_k)_c}.
\end{align}
\end{enumerate}
\end{proposition}

Before proving Proposition~\ref{pr:universal-ca-c}, let us show
that it implies the desired statements \eqref{eq:cluster-vars-principal-list}
and \eqref{eq:primitive-exchange-principal}.
Indeed, let $\PP^\circ = \Trop(y_i^\circ: i \in I)$ be the coefficient semifield
of the cluster algebra with principal coefficients
$\Acal^\circ(c)$ (see Definition~\ref{def:principal-coeffs}).
Let $\varphi: \widetilde \PP(c) \to \PP^\circ$ be the tropical
semifield homomorphism acting on the generators as follows:
\begin{equation}
\label{eq:principal-specialization-of universal-coeffs}
\varphi(p[\gamma]) = \begin{cases}
y_i & \text{if $\gamma = \omega_i$;} \\[.05in]
1 & \text{otherwise.}
\end{cases}
\end{equation}
Applying Proposition~\ref{pr:coef-specialization} to this
homomorphism, we see that the cluster variables in $\Acal^\circ(c)$
are in a natural bijection with those in $\widetilde \Acal(c)$,
and so can be labeled by the same set $\Pi(c)$.
If we denote by $x_{\gamma;c}$ the cluster variable in
$\Acal^\circ(c)$ corresponding to $x[\gamma] \in \widetilde
\Acal(c)$, then in view of
Proposition~\ref{pr:coef-specialization},
both \eqref{eq:cluster-vars-principal-list}
and \eqref{eq:primitive-exchange-principal} become consequences of
the following statement: under the homomorphism~$\varphi$, the
relations \eqref{eq:exchange-primitive-special-universal}
and \eqref{eq:exchange-primitive-nonspecial-universal} specialize
respectively to \eqref{eq:exchange-primitive-special}
and \eqref{eq:exchange-primitive-nonspecial}.
But this is immediate from the definition of the
$c$-compatibility degree.

\begin{proof}[Proof of Proposition~\ref{pr:universal-ca-c}]
We obtain the required assertion as a consequence of the following
two statements:
\begin{enumerate}
\item For the bipartite Coxeter element $t$, the cluster algebra $\widetilde \Acal(t)$
with required properties can be identified (by renaming the
cluster variables and the generators of the coefficient semifield)
with the cluster algebra with \emph{universal coefficients}
constructed in \cite[Theorem~12.4]{ca4}.
\item If $\widetilde \Acal(c)$ is the cluster algebra with required properties
for $c = s_1 \cdots s_n$, and $\tilde c = s_2 \cdots s_n s_1$, then
$\widetilde \Acal(\tilde c)$ can be identified (as above) with $\widetilde \Acal(c)$.
\end{enumerate}

Thus, each $\widetilde \Acal(c)$ for different choices of a Coxeter element~$c$
is a realization of the same cluster algebra with universal
coefficients, but with a different choice of an acyclic initial seed (with
the exchange matrix $B(c)$).

To prove (1), let us first recall the nomenclature of cluster
variables in \cite{ca4}.
Following \cite[Definition~9.1]{ca4}, for every $i \in I$ and $m \in \ZZ$ such that
$\varepsilon(i) = (-1)^m$, we define a root $\alpha(i;m)$, by setting, for all $r \geq 0$:
\begin{align}
\label{eq:dems-m-positive}
 \alpha(i;r) &=
\underbrace{t_- t_+ \cdots t_{\varepsilon(i)}}_{r
  \text{~factors}}(-\alpha_i) \qquad\, \text{for $\varepsilon(i) = (-1)^r$;}\\
\label{eq:dems-m-negative}
 \alpha(j;-r-1) &=
\underbrace{t_+ t_- \cdots t_{\varepsilon(j)}}_{r
\text{~factors}}(-\alpha_j) \qquad \text{for $\varepsilon(j) = (-1)^{r-1}$.}
\end{align}
In particular, we have
$\alpha(i;m) = - \alpha_i$ for $m \in \{0, -1\}$.

The following proposition is a consequence of a classical
result of R.~Steinberg \cite{steinberg} (cf.~\cite[Lemma~2.1,
Proposition~2.5]{yga}).
\footnote{Unfortunately, the corresponding result in
\cite[Proposition~9.3]{ca4} was stated incorrectly.}

\begin{proposition}
\label{pr:alpha-i-m}
The roots $\alpha(i;m)$ with $m \in [-h-1,-2]$, and $\varepsilon(i) = (-1)^m$
are positive, distinct, and every positive root is of such a form.
Furthermore, for $m \in [0,h+1]$, and $\varepsilon(i) = (-1)^m$,
we have $\alpha(i;m) = \alpha(i^\star;m-h-2)$.
\end{proposition}

In \cite{ca4}, the cluster variables are denoted by $x_{i;m}$ with
$i \in I$, $m \in \ZZ$ and $\varepsilon(i) = (-1)^m$.
If we parameterize them by the set $\Phi_{\geq -1}$ as in
Section~\ref{sec:Cluster-algebras} above (that is, using denominator vectors),
then by \cite[Theorem~10.3]{ca4},
each $x_{i;m}$ for $m \in [-h-1, h+1]$ will correspond to the root
$\alpha(i;m)$, and so we will write $x_{i;m} = x[\alpha(i;m)]$.
Comparing \eqref{eq:dems-m-negative} with the description of the
bijection $\psi: \Pi(t) \to \Phi_{\geq -1}$ given in
Lemma~\ref{lem:bijection-almost-positive-roots-Pit}, it is easy to see that,
for $0 \leq m \leq h(i;t)$, we have
\begin{equation}
\label{eq:psi-i-m}
\psi(t^m \omega_i) = \begin{cases}
\alpha(i;-2m) & \text{if $\varepsilon(i) = +1$;} \\[.05in]
\alpha(i;-2m-1) & \text{if $\varepsilon(i) = -1$.}
\end{cases}
\end{equation}

According to \cite[(8.12), (10.11)]{ca4}, the primitive exchange
relations are those having the product $x_{j;m-1}
x_{j;m+1}$ in the left hand side.
Renaming each cluster variable $x_{i;m}$ in these relations as
$x[\gamma]$, where $\gamma \in \Pi(t)$ is such that $\psi(\gamma)
= \alpha(i;m)$, and ignoring the coefficients, an easy check using \eqref{eq:psi-i-m}
shows that the cluster variables appear in these relations in exactly the same way as
in \eqref{eq:exchange-primitive-special-universal} and \eqref{eq:exchange-primitive-nonspecial-universal}.

To deal with coefficients, recall that the generators of the coefficient
semifield in \cite[Theorem~12.4]{ca4} are of the form
$p[\alpha^\vee]$ with $\alpha^\vee \in \Phi_{\geq - 1}^\vee$.
We identify $\Phi_{\geq - 1}^\vee$ with $\Pi(t)$ by means of the
following modification of Lemma~\ref{lem:bijection-almost-positive-roots-Pit}:

\begin{lemma}
\label{lem:bijection-almost-positive-coroots-Pit}
There is a bijection $\psi^\vee: \Pi(t) \to \Phi_{\geq -1}^\vee$ uniquely
determined by the properties that $\psi^\vee(\omega_i) = \varepsilon(i) \alpha_i^\vee$
for all $i \in I$, and $\psi^\vee \circ \tau_t = \tau^{-1} \circ \psi^\vee$.
Explicitly, $\psi^\vee$ is given by
\begin{equation}
\label{eq:psi-vee}
\psi^\vee(\gamma) = (\tau_+ \psi(\gamma))^\vee,
\end{equation}
where $\psi$ is as in
Lemma~\ref{lem:bijection-almost-positive-roots-Pit}.
\end{lemma}

\begin{proof}
The uniqueness statement is clear, as well as the fact that
\eqref{eq:psi-vee} defines a bijection $\Pi(t) \to \Phi_{\geq -1}^\vee$.
Using Lemma~\ref{lem:bijection-almost-positive-roots-Pit} and
\eqref{eq:tau-action}, we obtain that this bijection satisfies
$$\psi^\vee(\omega_i) = (\tau_+ (- \alpha_i))^\vee =
\varepsilon(i) \alpha_i^\vee,$$
and
$$\psi^\vee \tau_t(\gamma) = (\tau_+ \tau \psi(\gamma))^\vee
= (\tau_- \psi(\gamma))^\vee = (\tau^{-1} \tau_+ \psi(\gamma))^\vee
= \tau^{-1} \psi^\vee(\gamma),$$
as desired.
\end{proof}

As a consequence of \eqref{eq:t-comp}, \eqref{eq:compatibility-2}
and \eqref{eq:psi-vee}, we have
\begin{equation}
\label{eq:psi-vee-t-comp}
(\gamma \| \delta)_t = (\psi^\vee(\delta) \| \psi^\vee(\gamma))
\quad (\gamma, \delta \in \Pi(t)),
\end{equation}

Now everything is ready to check the following: replacing
$p[\gamma]$ with $p[\psi^\vee(\gamma)]$ for all $\gamma \in \Pi(t)$ transforms each element
$y_j$ given by \eqref{eq:y-initial-c-universal} (with $c=t$)
into the element
$$y_{j;0} = \prod_{\alpha^\vee \in \Phi^\vee_{\geq -1}}
p[\alpha^\vee]^{\varepsilon(j) [\alpha^\vee : \alpha_j^\vee]}$$
in \cite[(12.5)]{ca4}.

\smallskip

{\bf Case 1.} Let $\varepsilon (j) = +1$.
Then we have $\psi^\vee(\omega_j) = \alpha_j^\vee$ and
$\psi^\vee(t\omega_j) = -\alpha_j^\vee$.
Since there are no indices $i$ with $i \prec_t j$,
if $\psi^\vee (\gamma) = \alpha^\vee \neq \pm \alpha_j^\vee$, then the
exponent of $p[\gamma]$ in the right hand side of \eqref{eq:y-initial-c-universal}
is equal to
$$(\gamma \| t \omega_j)_t = (- \alpha_j^\vee \| \alpha^\vee) =
[\alpha^\vee : \alpha_j^\vee].$$
Thus, the replacement of $p[\gamma]$ with $p[\alpha^\vee]$
indeed transforms $y_j$ into $y_{j;0}$.

\smallskip

{\bf Case 2.} Let $\varepsilon (j) = -1$.
Then we have $\psi^\vee(\omega_j) = -\alpha_j^\vee$ and
$\psi^\vee(t\omega_j) = \alpha_j^\vee$.
In this case, there are no indices $i$ with $j \prec_t i$, so by
Lemma~\ref{lem:comp-c-identity}, if $\psi^\vee (\gamma) = \alpha^\vee \neq \pm \alpha_j^\vee$,
then the exponent of $p[\gamma]$ in the right hand side of \eqref{eq:y-initial-c-universal}
is equal to
$$-(\gamma \| \omega_j)_t = - (- \alpha_j^\vee \| \alpha^\vee) =
- [\alpha^\vee : \alpha_j^\vee],$$
proving our claim in this case as well.

\smallskip

To finish the proof of (1), it remains to check that the same
replacement of each $p[\gamma]$ with $p[\psi^\vee(\gamma)]$
transforms the coefficients in the relations
\eqref{eq:exchange-primitive-special-universal} and
\eqref{eq:exchange-primitive-nonspecial-universal} into
the coefficients of the corresponding relations in the cluster
algebra from \cite[Theorem~12.4]{ca4}.
Pick some $\delta \in \Pi(t)$, and let
$\alpha^\vee = \psi^\vee(\delta)$.
Recall from \cite[Section~12]{ca4} that $p[\alpha^\vee]$ is the
primitive coefficient (the one at the non-trivial product of
cluster variables) in the relation between $x[\tau_- \alpha]$ and
$x[\tau_+ \alpha]$.
This agrees with the fact that $p[\delta]$ is the primitive coefficient in the relation
(of the form \eqref{eq:exchange-primitive-special-universal} or
\eqref{eq:exchange-primitive-nonspecial-universal}) between $x[\tau_t^{-1}\delta]$ and $x[\delta]$.
It remains to observe that the constant term in the same relation
\eqref{eq:exchange-primitive-special-universal} or
\eqref{eq:exchange-primitive-nonspecial-universal}, that is, the
element $\prod_{\gamma \in \Pi(t)} p[\gamma]^{(\gamma \| \delta)_t}$,
transforms, in view of \eqref{eq:psi-vee-t-comp}, into
$\prod_{\beta^\vee \in {\Phi_{\geq -1}^\vee}}
p[\beta^\vee]^{(\alpha^\vee \| \beta^\vee)}$,
which is the constant term in the
corresponding relation in \cite{ca4} (see \cite[(12.18)]{ca4},
where unfortunately there is an annoying typo: the correct exponent on
the right is $(\alpha^\vee \| \beta^\vee)$ instead of $(\beta^\vee \|
\alpha^\vee)$).
This completes the check of property (1) above.

\smallskip

To prove (2), it is enough to show the following: after replacing,
for each $\gamma \in \Pi(\tilde c)$, the element $x[\gamma] \in
\widetilde \Acal(\tilde c)$ by $x[\psi_{c, \tilde c}(\gamma)] \in
\widetilde \Acal(c)$, and the coefficient $p[\gamma] \in
\widetilde \PP(\tilde c)$ by $p[\psi_{c, \tilde c}(\gamma)] \in
\widetilde \PP(c)$ (see Lemma~\ref{lem:bijection-Pic-Pitildec}),
the properties (3) and (4) for $\widetilde \Acal(\tilde c)$
in Proposition~\ref{pr:universal-ca-c}
become identical to the corresponding properties for $\widetilde \Acal(c)$.

Starting with property (4), it is enough to check it for the
relations \eqref{eq:exchange-primitive-special-universal}, since
the relations \eqref{eq:exchange-primitive-nonspecial-universal} are
obtained from them by repeatedly applying the transformation
$\tau_{\tilde c}$ (resp. $\tau_{c}$) to all occurring labeling weights $\gamma \in
\Pi(\tilde c)$ (resp. $\gamma \in \Pi(c)$), and the bijection
$\psi_{c, \tilde c}: \Pi(\tilde c) \to \Pi(c)$ intertwines $\tau_{\tilde c}$ with $\tau_{c}$.
The check that every relation \eqref{eq:exchange-primitive-special-universal}
for $\widetilde \Acal(\tilde c)$ turns into one of the \eqref{eq:exchange-primitive-special-universal}
or \eqref{eq:exchange-primitive-nonspecial-universal} for $\widetilde \Acal(c)$
after applying $\psi_{c, \tilde c}$ is immediate from the definitions.

To deal with the remaining property (3), we note that $\psi_{c, \tilde c}$ sends
the initial cluster $(x[\omega_1], \dots, x[\omega_n])$ in $\widetilde \Acal(\tilde c)$
into the cluster $(x[c\omega_1], x[\omega_2], \dots,
x[\omega_n])$ in $\widetilde \Acal(c)$ obtained from the initial
one by the mutation $\mu_1$.
Remembering \eqref{eq:Bc} and \eqref{eq:matrix-mutation}, we see
that $\mu_1$ sends the exchange matrix $B(c)$ into $B(\tilde c)$.
It remains to show that the transformation \eqref{eq:y-mutation}
(for $k=1$) applied to the coefficient tuple in $\widetilde
\PP(c)$ given by \eqref{eq:y-initial-c-universal} produces the
elements $y'_j$ obtained from the corresponding $\tilde y_j \in
\widetilde \PP(\tilde c)$ by applying $\psi_{c, \tilde c}$ to all
the labeling weights.

For $j=1$, we can express $y_1$ as
$$y_1 =  p[c \omega_1]^{-1} \prod_{\gamma \in \Pi(c)}
p[\gamma]^{(\gamma \| c \omega_1)_c};$$
therefore, we have
$$y'_1 = y_1^{-1} = p[c \omega_1]
\prod_{\gamma \in \Pi(c)}
p[\gamma]^{-(\gamma \| c \omega_1)_c}.$$
On the other hand, using \eqref{eq:comp-c-omega-omega}, we can
express $\tilde y_1 \in \widetilde \PP(\tilde c)$ as
$$\tilde y_1 = p[\omega_1] \prod_{\gamma \in \Pi(\tilde c)}
p[\gamma]^{-(\gamma \| \omega_1)_{\tilde c}},$$
which indeed transforms into $y'_1$ by applying $\psi_{c, \tilde c}$ to all
its labeling weights.

Finally, for $j > 1$, we have $b_{1,j} = - a_{1,j} \geq 0$, hence
\eqref{eq:y-mutation} results in
\begin{align*}
y'_j &= y_j y_1^{-a_{1,j}}  p[c \omega_1]^{-a_{1,j}}
= y_j \prod_{\gamma \in \Pi(c)} p[\gamma]^{- a_{1,j} (\gamma \| c
\omega_1)_c}\\
&= p[\omega_j]  p[c \omega_j]^{-1}
\prod_{\gamma \in \Pi(c) - \{\omega_j, c \omega_j\}}
p[\gamma]^{(\gamma \| c \omega_j)_c + \sum_{1 \neq i \prec_c j} a_{i,j} (\gamma \| c
\omega_i)_c},
\end{align*}
which is again obtained from $\tilde y_j$ by applying $\psi_{c, \tilde c}$ to all
its labeling weights.

This finishes the proof of  Proposition~\ref{pr:universal-ca-c},
hence also the proofs of Theorems~\ref{th:ca-general-2}
and~\ref{th:ca-general-3}.
\end{proof}

In the rest of the section we give some corollaries of the
developed techniques.
First, since the $c$-compatibility degree for every Coxeter element~$c$
is obtained by a chain of bijections from the compatibility degree in
\cite[Section~3.1]{yga} (see \eqref{eq:t-comp} and \eqref{eq:tilde-c-thru-c-comp}),
combining the above results with
Proposition~\ref{pr:root-clusters} yields the following corollary.

\begin{corollary}
\label{cor:c-clusters}
For an arbitrary Coxeter element~$c$ and any $\gamma, \delta \in
\Pi(c)$, the cluster variables $x_{\gamma;c}$ and $x_{\delta;c}$
in $\Acal(c)$ belong to the same cluster if and only if
$(\gamma \| \delta)_c = 0$ (in particular, the latter condition is
symmetric in $\gamma$ and $\delta$).
\end{corollary}

Second, combining the arguments in the proof of Proposition~\ref{pr:universal-ca-c}
with Proposition~\ref{pr:principal-embedding}, we conclude that
all the algebras with principal coefficients $\Acal(c)$ for
different choices of the Coxeter element~$c$ are isomorphic to each other.
To be more specific, let us again suppose that $c = s_1 s_2 \cdots
s_n$ and $\tilde c = s_2 \cdots s_n s_1$.
Then the isomorphism $\widetilde \Acal(\tilde c) \to \widetilde
\Acal(c)$ constructed above gives rise to the isomorphism
$\Acal(\tilde c) \to \Acal(c)$ given as follows.

\begin{corollary}
\label{cor:c-tildec-isom-algebraic}
The correspondence in Proposition~\ref{pr:principal-embedding}
establishes an isomorphism of cluster algebras $\varphi: \Acal(\tilde c) \to \Acal(c)$
acting on the initial cluster variables $x_{\omega_i;\tilde c}$
and the initial coefficient tuple $(y_{1;\tilde c},
\dots, \tilde y_{n;\tilde c})$ in $\Acal(\tilde c)$ as follows:
$$\varphi (x_{\omega_i;\tilde c}) = x_{\omega_i;c}, \quad
\varphi (y_{i;\tilde c}) = y_{i;c}  y_{1;c}^{-a_{1,i}} \quad (i
\neq 1),$$
$$\varphi (x_{\omega_1;\tilde c}) = x_{c\omega_1;c} =
x_{\omega_1;c}^{-1} (y_{1;c} + \prod_{i \neq 1}
x_{\omega_i;c}^{-a_{i,1}}), \quad \varphi (y_{1;\tilde c}) =
y_{1;c}^{-1}.$$
\end{corollary}

\begin{remark}
\label{rem:c-tildec-isom-geometric}
Since, in view of Theorem~\ref{th:ca-general-1}, the algebra
$\Acal(c)$ is the coordinate ring of the variety $L^{c,c^{-1}}$,
the isomorphism $\varphi: \Acal(\tilde c) \to \Acal(c)$ gives rise
to a biregular isomorphism
$\varphi^*: L^{c,c^{-1}} \to L^{\tilde c,\tilde c^{-1}}$.
This isomorphism can be described as follows.
Note that the definition \eqref{eq:reduced cell} and the standard
properties of double Bruhat cells imply the following statement:
if in the Weyl group~$W$ we have factorizations $u = u_1 \cdots u_k$ and $v =
v_1 \cdots v_k$ such that $\ell(u) = \ell(u_1) + \cdots +
\ell(u_k)$ and $\ell(v) = \ell(v_1) + \cdots + \ell(v_k)$, then
the product map in~$G$ induces an open embedding
$$L^{u_1, v_1} \times \cdots \times L^{u_k, v_k} \hookrightarrow L^{u, v}.$$
In particular, denoting $c_\circ = s_2 \cdots s_n$, we have open embeddings
$$L^{s_1,e} \times L^{c_\circ, c_\circ^{-1}} \times L^{e,s_1} \hookrightarrow
L^{c,c^{-1}}, \quad L^{e,s_1} \times L^{c_\circ, c_\circ^{-1}} \times L^{s_1,e} \hookrightarrow
L^{\tilde c,\tilde c^{-1}}.$$
Note also that
$$L^{s_1,e} = \{x_{-1}(u): u \in \CC^*\}, \quad L^{e,s_1} = \{x_{1}(t): t \in \CC^*\}\ .$$
Now we claim that~$\varphi^*$ restricts to an isomorphism
$$L^{s_1,e} \cdot L^{c_\circ, c_\circ^{-1}} \cdot L^{e,s_1} \to
L^{e,s_1} \cdot L^{c_\circ, c_\circ^{-1}} \cdot L^{s_1,e}$$
given by
\begin{equation}
\label{eq:c-tildec-isom-geometric}
\varphi^*(x_{-1}(u) x_\circ x_1(t)) = x_1(u) x_\circ
x_{-1}(t^{-1}) \quad (u, t \in \CC^*, \,\, x_\circ \in L^{c_\circ,
c_\circ^{-1}}) \ .
\end{equation}
It suffices to prove \eqref{eq:c-tildec-isom-geometric} for $x_\circ$ of the form
$x_{-2}(u_2) \cdots x_{-n}(u_n) x_{n}(t_n)
\cdots x_{2}(t_2)$ with all $u_i$ and $t_i$ nonzero complex numbers, since such elements
form an open dense subset of $L^{c_\circ, c_\circ^{-1}}$.
Then \eqref{eq:c-tildec-isom-geometric} can be checked by a straightforward calculation using
\eqref{eq:initial-principal-minors}, \eqref{eq:initial-yj}
and the commutation relations in \cite[Proposition~7.2]{bz01}; we
leave the details to the reader.
\end{remark}

\section{Proofs of Theorems~\ref{th:denom-vector}, \ref{th:g-vector} and \ref{th:F-poly},
and their corollaries}
\label{sec:proofs-g-vector-F-poly}

In this section we work with the cluster algebra $\Acal(c)$ from
Theorem~\ref{th:ca-general-1}.
When it is convenient, we assume without further warning
that the index set $I$ is $[1,n]$, and the Coxeter
element~$c$ under consideration is $s_1 \cdots s_n$.

\begin{proof}[Proof of Theorem~\ref{th:denom-vector}]
For $\gamma \in \Pi(c)$, let $\dd(\gamma)$ denote the denominator
vector of the cluster variable $x_{\gamma;c}$ with respect to the
initial cluster $(x_{\omega_i;c} : i \in I)$.
We identify $\ZZ^n$ with the root lattice~$Q$ using the basis of simple roots,
and so assume that $\dd(\gamma) \in Q$.
In view of \cite[(7.6), (7.7)]{ca4}, the vectors
$\dd(\gamma)$ are uniquely  determined from the initial conditions
\begin{equation}
\label{eq:denoms-initial}
\dd(\omega_i) = - \alpha_i,
\end{equation}
and the recurrence relations
\begin{equation}
\label{eq:denoms-recurrence}
\dd(c^m \omega_j) + \dd(c^{m-1} \omega_j) = [- \sum_{i \prec_c j}
a_{i,j} \dd(c^m \omega_i) - \sum_{j \prec_c i}
a_{i,j} \dd(c^{m-1} \omega_i)]_+
\end{equation}
(for all $j \in I$ and $1 \leq m \leq h(j;c)$),
which follow from \eqref{eq:exchange-primitive-nonspecial}; here
the notation $[v]_+$ for $v \in Q$ is understood component-wise, i.e.,
$$[\sum a_i \alpha_i]_+ = \sum [a_i]_+ \alpha_i.$$
We need to show that the solution of the relations \eqref{eq:denoms-recurrence} with
the initial conditions \eqref{eq:denoms-initial} is given by
\begin{equation}
\label{eq:denoms-answer}
\dd(c^m \omega_j) = c^{m-1}\omega_j - c^m \omega_j = c^{m-1} \beta_j  \quad (j \in I, 1 \leq m \leq h(j;c))
\end{equation}
(see Lemma~\ref{lem:two-identities}).

First let us check that the values given by
\eqref{eq:denoms-initial} and \eqref{eq:denoms-answer} satisfy
\eqref{eq:denoms-recurrence} for $m=1$.
Indeed, the right hand side of \eqref{eq:denoms-recurrence} is
equal to
$$[- \sum_{i \prec_c j} a_{i,j} \beta_i + \sum_{j \prec_c i}
a_{i,j} \alpha_i]_+ = - \sum_{i \prec_c j} a_{i,j} \beta_i,$$
since each $\beta_i = s_1 \cdots s_{i-1} \alpha_i$ is a positive
linear combination of the roots $\alpha_{i'}$ with $i' < i$
(in any total order on the index set~$I$ compatible with the
relation $i' \prec_c i$).
On the other hand, the left hand side of \eqref{eq:denoms-recurrence} is
equal to $\beta_j - \alpha_j$.
So the two sides are equal to each other by
\eqref{eq:beta-telescoping}.

For $m \geq 2$, the desired identity
\eqref{eq:denoms-recurrence} takes the form
$$c^{m-1} \beta_j + c^{m-2} \beta_j =
- \sum_{i \prec_c j}
a_{i,j}c^{m-1} \beta_i  - \sum_{j \prec_c i}
a_{i,j} c^{m-2} \beta_i,$$
which can be simplified to
\begin{equation}
\label{eq:intermediate-identity}
\beta_j + \sum_{i \prec_c j} a_{i,j} \beta_i =
- c^{-1}(\beta_j + \sum_{j \prec_c i} a_{i,j} \beta_i).
\end{equation}
Now the left hand side of \eqref{eq:intermediate-identity}
is equal to $\alpha_j$ by \eqref{eq:beta-telescoping}.
On the other hand, since
$$- c^{-1} \beta_i = - c^{-1} s_1 \cdots s_{i-1} \alpha_i =
- s_n \cdots s_{i+1} s_i \alpha_i = s_n \cdots s_{i+1} \alpha_i,$$
the right hand side of \eqref{eq:intermediate-identity} is also
equal to $\alpha_j$ by the same equality \eqref{eq:beta-telescoping}
with the Coxeter element $c$ replaced by $c^{-1}$.
This completes the proof of Theorem~\ref{th:denom-vector}.
\end{proof}

\begin{proof}[Proof of Corollary~\ref{cor:denoms}]
The first assertion is clear, and the second one follows at once by combining the expression for
denominator vectors given in Theorem~\ref{th:denom-vector} with
Corollary~\ref{cor:c-clusters} and
formula (2) in Proposition~\ref{pr:comp-degree-Pic}.
\end{proof}

\begin{proof}[Proof of Theorem~\ref{th:g-vector}]
Recall that the $\gg$-vector $\gg_{z;\xx} = (g_1, \dots, g_n) \in \ZZ^n$
and the $F$-polynomial $F_{z;\xx}(t_1, \dots, t_n) \in \ZZ[t_1, \dots, t_n]$
of a  cluster variable~$z \in \Acal^\circ$ with respect to the initial cluster~$\xx$
are defined by \eqref{eq:cluster-vars-principal}.
As prescribed by Theorem~\ref{th:g-vector}, we now identify
$\ZZ^n$ with the weight lattice~$P$ using the basis of fundamental weights.
Following \cite{ca4}, we introduce the $P$-(multi)grading in the
Laurent polynomial ring $\ZZ[\xx^{\pm 1}, \yy^{\pm 1}]$
by setting
\begin{equation}
\label{eq:deg-x-y}
{\rm deg} (x_{\omega_j;c}) = \omega_j, \quad {\rm deg} (y_{j;c}) =
- \sum_{i \in I} b_{i,j} \omega_i = \sum_{i \prec_c j} a_{i,j} \omega_i
- \sum_{j \prec_c i} a_{i,j} \omega_i.
\end{equation}
Then by \eqref{eq:y-hat}, all elements $\widehat y_j^\circ$ are
homogeneous of degree~$0$, hence the same is true for
the factor $F_{z;\xx^\circ}(\widehat y_1^\circ, \dots, \widehat y_n^\circ)$
in \eqref{eq:cluster-vars-principal}, so the $\gg$-vector is equal to
$\gg_{z;\xx} = {\rm deg}(z)$.

Following \eqref{eq:Lc-c-1-factor}, let us write a generic element of $L^{c,c^{-1}}$ in
the form
$$x_{-1}(u_1) \cdots x_{-n}(u_n) x_n(t_n) \cdots x_1(t_1),$$
and view all $t_j$ and
$u_j$ as rational functions on $L^{c,c^{-1}}$.
In view of \eqref{eq:initial-principal-minors} and
\eqref{eq:initial-yj}, the functions $t_j$ and $u_j$ are Laurent
monomials in the initial cluster variables $x_{\omega_i;c}$ and
the coefficients $y_{i;c}$, and we have
\begin{align}
\label{eq:deg-u-t}
&{\rm deg} (u_j) = -{\rm deg} (x_{\omega_j;c}) = -\omega_j, \\
\nonumber
&{\rm deg} (t_j) = {\rm deg} (y_{j;c}) + {\rm deg} (u_j) +
\sum_{i \prec_c j} a_{i,j}  {\rm deg} (u_i) =
- (\omega_j + \sum_{j \prec_c i} a_{i,j} \omega_i).
\end{align}

Remembering \eqref{eq:x-gamma-c}, Theorem~\ref{th:g-vector}
can be restated as follows: with respect to the $P$-grading
given by \eqref{eq:deg-u-t}, we have
\begin{equation}
\label{eq:deg-Delta}
{\rm deg}(\Delta_{\gamma, \gamma}
(x_{-1}(u_1) \cdots x_{-n}(u_n) x_n(t_n) \cdots x_1(t_1))) = \gamma.
\end{equation}
To prove \eqref{eq:deg-Delta}, we write each factor $x_{-j}(u_j)$ as
$x_{\bar j}(u_j) u_j^{- \alpha_j^\vee}$ (see \eqref{eq:xnegative}), and move all the
factors $u_j^{- \alpha_j^\vee}$ all the way to the right, using
the commutation relations in (\ref{eq:H-conjugation}).
We obtain
$$x_{-1}(u_1) \cdots x_{-n}(u_n) x_n(t_n) \cdots x_1(t_1) =
x_{\bar 1}(w_1) \cdots x_{\bar n}(w_n) x_n(v_n) \cdots x_1(v_1)
\prod_{i \in I} u_i^{- \alpha_i^\vee},$$
where the $w_j$ and $v_j$ are given by
\begin{equation}
\label{eq:w-v-thru-u-t}
w_j = u_j \prod_{i \prec_c j} u_i^{a_{i,j}}, \quad
v_j = t_j \prod_{i = 1}^n u_i^{-a_{i,j}} \ .
\end{equation}
Expressing $\gamma$ in the form $\gamma = \sum_i g_i \omega_i$ and
using \eqref{eq:minor-weight}, we obtain
\begin{align*}
&\Delta_{\gamma, \gamma} (x_{-1}(u_1) \cdots x_{-n}(u_n) x_n(t_n) \cdots
x_1(t_1)))\\
& = \Delta_{\gamma, \gamma}
(x_{\bar 1}(w_1) \cdots x_{\bar n}(w_n) x_n(v_n) \cdots x_1(v_1))
\cdot (\prod_i u_i^{- \alpha_i^\vee})^\gamma\\
& = \Delta_{\gamma, \gamma}
(x_{\bar 1}(w_1) \cdots x_{\bar n}(w_n) x_n(v_n) \cdots x_1(v_1))
\prod_i u_i^{- g_i},
\end{align*}
hence
\begin{align*}
&{\rm deg}(\Delta_{\gamma, \gamma}
(x_{-1}(u_1) \cdots x_{-n}(u_n) x_n(t_n) \cdots x_1(t_1)))\\
&= \gamma + {\rm deg}(\Delta_{\gamma, \gamma}
(x_{\bar 1}(w_1) \cdots x_{\bar n}(w_n) x_n(v_n) \cdots
x_1(v_1))).
\end{align*}
Thus it remains to prove that
$${\rm deg}(\Delta_{\gamma, \gamma}
(x_{\bar 1}(w_1) \cdots x_{\bar n}(w_n) x_n(v_n) \cdots x_1(v_1))) = 0 \ .$$

Using \eqref{eq:w-v-thru-u-t} and \eqref{eq:deg-u-t}, we deduce that
$${\rm deg} (w_j) = -\omega_j - \sum_{i \prec_c j} a_{i,j} \omega_i = -{\rm deg} (v_j) \ .$$
Now for every nonzero complex numbers $w_1, \dots, w_n$ there is
$a \in H$ such that $w_j = a^{\alpha_j}$ for all~$j$.
Using \eqref{eq:minor-weight} and \eqref{eq:H-conjugation}, we
conclude that
\begin{align*}
&\Delta_{\gamma, \gamma}
(x_{\bar 1}(w_1) \cdots x_{\bar n}(w_n) x_n(v_n) \cdots
x_1(v_1))\\
&= \Delta_{\gamma, \gamma}
(a x_{\bar 1}(w_1) \cdots x_{\bar n}(w_n) x_n(v_n) \cdots
x_1(v_1)a^{-1})\\
&= \Delta_{\gamma, \gamma}
(x_{\bar 1}(1) \cdots x_{\bar n}(1) x_n(w_n v_n) \cdots x_1(w_1
v_1)).
\end{align*}
Since each product $w_j v_j$ has degree~$0$,
this completes the proof of Theorem~\ref{th:g-vector}.
\end{proof}

\begin{proof}[Proof of Theorem~\ref{th:F-poly}]
We use the notation in the proof of Theorem~\ref{th:g-vector}.
In view of \eqref{eq:cluster-vars-principal}, the $F$-polynomial $F_{z;\xx}(t_1, \dots, t_n)$
is obtained from the expansion of the cluster variable $z = x_{\gamma;c}$
by specializing all initial cluster variables $x_{\omega_i;c}$
to~$1$, and each coefficient~$y_{i;c}$ to~$t_i$.
In view of \eqref{eq:initial-principal-minors} and
\eqref{eq:initial-yj}, this amounts to specializing all~$u_i$
to~$1$ in $\Delta_{\gamma, \gamma} (x_{-1}(u_1) \cdots x_{-n}(u_n) x_n(t_n) \cdots
x_1(t_1)))$, which is exactly the desired assertion.
\end{proof}

\begin{proof}[Proof of Corollary~\ref{cor:ct-1}]
The constant term of the polynomial $F_{z;\xx}(t_1, \dots, t_n)$
is obtained by specializing all~$t_i$
to~$0$ in \eqref{eq:F-poly}, i.e., is equal to $\Delta_{\gamma, \gamma}
(x_{\bar 1}(1) \cdots x_{\bar n}(1))$.
The fact that the latter is equal to~$1$ follows at once by the
representation-theoretic arguments in the above proof of \eqref{eq:cmcm+1inati}.
\end{proof}

\begin{remark}
\label{rem:cambrian-2}
We have already mentioned in Remark~\ref{rem:cambrian-1} that
the set $\Pi(c)$ is identical with the set of generators of the
$c$-Cambrian fan studied in \cite{reading-speyer}.
In particular, our Theorem~\ref{th:g-vector} implies
\cite[Theorem~10.2]{reading-speyer}, proved there only modulo \cite[Conjecture~7.12]{ca4}.
Comparing our Corollary~\ref{cor:c-clusters} with
\cite[Theorem~10.1]{reading-speyer}, we obtain the following
description of the cones in the $c$-Cambrian fan: they
are exactly the cones generated by subsets $\{\gamma_1, \dots,
\gamma_k\} \subset \Pi(c)$ such that $(\gamma_i \| \gamma_j)_c =
0$ for all $i,j = 1, \dots, k$.
\end{remark}

\section{Proof of Theorem~\ref{th:ca-principal-minors-A-special}}
\label{sec:proof-of-ca-principal-minors-A-special}

In this section we assume that $G = SL_{n+1}(\CC)$ is of type $A_n$, and the
Coxeter element $c$ is equal to $s_1 \cdots s_n$ in the standard numbering of simple roots.
As usual, the Weyl group $W$ is identified with the symmetric
group~$S_{n+1}$, so that $s_i$ becomes a simple
transposition~$(i,i+1)$.
Then~$c$ is the cycle $(1,2, \cdots, n+1)$.

In this case, the generalized minors $\Delta_{\gamma,\delta}$ specialize to the ordinary minors
(i.e., determinants of square submatrices) as follows.
The weights in $W \omega_k$ are in bijection with the $k$-subsets of $[1,n+1] = \{1,\dots,n+1\}$,
so that~$W = S_{n+1}$ acts on them in a natural way, and $\omega_k$ corresponds to~$[1,k]$.
If~$\gamma$ and~$\delta$ correspond to $k$-subsets $I$ and $J$, respectively, then
$\Delta_{\gamma,\delta} = \Delta_{I,J}$ is the minor with
the row set~$I$ and the column set~$J$.

Now let~$L$ denote the subvariety of~$G$ consisting of tridiagonal
matrices~$M$ of the form \eqref{eq:tridiagonal} (with all $y_1, \dots, y_n$ non-zero).
We start our proof of Theorem~\ref{th:ca-principal-minors-A-special} by showing
that $L$ is indeed the reduced double cell~$L^{c,c^{-1}}$.
We use the characterization of $L^{c,c^{-1}}$ given by Propositions~\ref{pr:double-cell-by-eqs}
and~\ref{pr:Luv-equations}.
It is well-known that the Bruhat order on $k$-elements subsets of
$[1,n+1]$ is component-wise: if $I = \{i_1 < \cdots < i_k\}$ and
$J = \{j_1 < \cdots < j_k\}$ then $I \leq J$ means that
$i_\nu \leq j_\nu$ for $\nu = 1, \dots, k$.
Thus, we can rewrite Propositions~\ref{pr:double-cell-by-eqs}
and~\ref{pr:Luv-equations} as follows.

\begin{corollary}
\label{cor:RDBC-SL}
A matrix $M \in SL_{n+1}(\CC)$ belongs to the reduced double Bruhat cell
$L^{c,c^{-1}}$ for $c = s_1 \cdots s_n$ if and only if
it satisfies the following conditions:
\begin{enumerate}
\item $\Delta_{I,[1,k]} = \Delta_{[1,k],I} = 0$
for $k = 1, \dots, n$, and all subsets
$I = \{i_1 < \cdots < i_k\} \subset [1,n+1]$
such that $i_\nu > \nu + 1$ for some $\nu = 1, \dots, k$.

\item $\Delta_{[1,k],[2,k+1]} \neq 0$, and
$\Delta_{[2,k+1],[1,k]} = 1$ for $k = 1, \dots, n$.
\end{enumerate}
\end{corollary}

To prove that $L = L^{c,c^{-1}}$, we start with the inclusion $L \subseteq L^{c,c^{-1}}$.
Let $M = (m_{i,j}) \in L$.
Since $M$ is tridiagonal, it satisfies condition (1) in
Corollary~\ref{cor:RDBC-SL} because every term
in the expansion of $\Delta_{I,[1,k]}(M)$  or
$\Delta_{[1,k],I}(M)$ (with~$I$ as in this condition) contains
at least one matrix entry $m_{i,j}$ with $|i-j| > 1$.
To prove condition (2), note that any tridiagonal matrix~$M$ becomes
triangular after removing the first column and the last row (or
the first row and the last column), implying that
\begin{align}
\label{eq:tridiag-product-minors}
&\Delta_{[1,k],[2,k+1]}(M) = m_{1,2} m_{2,3} \cdots m_{k,k+1}\\
\nonumber
&\Delta_{[2,k+1],[1,k]}(M) = m_{2,1} m_{3,2} \cdots m_{k+1,k}.
\end{align}
So for $M \in L$, we have
$\Delta_{[1,k],[2,k+1]}(M) = y_1 \cdots y_k \neq 0$, and
$\Delta_{[2,k+1],[1,k]}(M) = 1$, as required.

To prove the reverse inclusion $L^{c,c^{-1}} \subseteq L$, first
let us show that every matrix in $L^{c,c^{-1}}$ is tridiagonal.
It is enough to check this for a generic element of the form
\eqref{eq:Lc-c-1-factor}.
Note that in our situation the matrix $x_{-k}(u)$ (resp. $x_k(t)$)
is obtained from the identity matrix by replacing the $2 \times 2$
submatrix with rows and columns $k$ and $k+1$ by
$\mat{u^{-1}}{0}{1}{u}$ (resp. $\mat{1}{t}{0}{1}$).
Then one can check that the product
$$x_{-1}(u_1) \cdots x_{-n}(u_n) x_n(t_n) \cdots x_1(t_1)$$
is tridiagonal by a direct matrix computation (or better yet, by
using the graphical formalism for computing determinants of such
products developed in \cite[Proof of Theorem~12]{fz-intel}).

Once we know that every matrix $M \in L^{c,c^{-1}}$ is tridiagonal,
the remaining conditions $m_{k,k+1} \neq 0$ and $m_{k+1,k} = 1$
follow at once from \eqref{eq:tridiag-product-minors} and
condition (2) in Corollary~\ref{cor:RDBC-SL}.
This concludes the proof of the equality  $L = L^{c,c^{-1}}$.

\smallskip

We continue the proof of Theorem~\ref{th:ca-principal-minors-A-special}.
Part (1) of Theorem~\ref{th:ca-principal-minors-A-special}
is a special case of Theorem~\ref{th:ca-general-1}.
We need only to observe that the expression for $y_{j;c}$ in \eqref{eq:y-j-c}
specializes to the restriction to~$L$ of the function
$\Delta_{[1,j],[2,j+1]}/\Delta_{[1,j-1],[2,j]}$, which, in view of \eqref{eq:tridiag-product-minors}
is equal to the above-diagonal matrix entry $y_j$ (see
\eqref{eq:tridiagonal}).

\smallskip

Part (2) of Theorem~\ref{th:ca-principal-minors-A-special}
is a special case of Theorem~\ref{th:ca-general-2}.
We need only to observe that in our case $h(k;c) = n+1-k$, and, for
$0 \leq m \leq n+1-k$,  the
weight $c^m \omega_k \in \Pi(c)$ corresponds to the subset
$[m+1, m+k] \subset [1,n+1]$.

Turning to Part (3), we first note that
the primitive exchange relations~\eqref{eq:exchange-primitive-special} and
\eqref{eq:exchange-primitive-nonspecial}
specialize to the following relations, which are among those in~\eqref{eq:exchange-A-principal}:
\begin{equation}
\label{eq:exchange-primitive-special-An}
x_{[1,k]} x_{[k+1, n+1]} =
y_{k} x_{[1,k-1]} x_{[k+2, n+1]} + 1;
\end{equation}
\begin{equation}
\label{eq:exchange-primitive-nonspecial-An}
x_{[m, m+k-1]}  x_{[m+1, m+k]} = y_{m}y_{m+1} \dots y_{m+k-1} + x_{[m, m+k]} x_{[m+1, m+k-1]},
\end{equation}
where $k \in [1,n]$, $1 \leq m \leq h(k;c) = n + 1 - k$.
To prove Theorem~\ref{th:ca-principal-minors-A-special}, it remains to show that
the exchange relations in $\Acal(c)$ are exactly those
in~\eqref{eq:exchange-A-principal}.

We start by recalling the geometric interpretation of the cluster
variables and clusters in type $A_n$ given in \cite[Section~12.2]{ca2}.
Namely, the cluster variables can be associated with the diagonals
of a regular $(n+3)$-gon $\Poly_{n+3}$, and the clusters are
all maximal collections of mutually non-crossing diagonals, which
are naturally identified with the \emph{triangulations} of $\Poly_{n+3}$.
We label the vertices of $\Poly_{n+3}$ by numbers $1, \ldots, n+3$ in the counter-clockwise
order, and denote by $x_{\langle i,j \rangle}$ the cluster variable associated
with the diagonal $\langle i,j \rangle$ connecting vertices $i$ and $j$
(with the convention that $x_{\langle i,j \rangle} = 1$ if $i$ and $j$ are two adjacent vertices
of $\Poly_{n+3}$).
According to \cite[(12.2)]{ca2}, in this realization the exchange relations have the form
\begin{equation}
\label{eq:ex-A-UC}
x_{\langle i,k \rangle}\, x_{\langle j, \ell \rangle}
 = p^+_{ik,j\ell}\, x_{\langle i,j \rangle}\, x_{\langle k, \ell \rangle}
+ p^-_{ik,j\ell}\, x_{\langle i,\ell \rangle}\, x_{\langle j, k \rangle} \ ,
\end{equation}
where $i,j,k,\ell$ are any four vertices of $\Poly_{n+3}$ taken in
counter-clockwise order, and $p^\pm_{ik,j\ell}$
are some elements of the coefficient semifield.

\begin{figure}[ht]
\setlength{\unitlength}{1.5pt}
\begin{picture}(160,120)(0,0)

\put(132,26){$1$}
\put(132,79){$2$}
\put(83,110){$3$}
\put(34,79){$4$}
\put(34,26){$5$}
\put(83,-3){$6$}

\put(114,59){$x_{[1,1]}$}
\put(90,54){$x_{[1,2]}$}
\put(57,33){$x_{[1,3]}$}
\put(108,84){$\textcolor[rgb]{1.00,0.00,0.00}{y_1}$}
\put(70,84){$\textcolor[rgb]{1.00,0.00,0.00}{y_2}$}
\put(48,54){$\textcolor[rgb]{1.00,0.00,0.00}{y_3}$}

\drawline(85,4)(129,30)(129,81.1077)(85,107.1077)(41,81.1077)(41,30)(85,4)

\light{
\drawline(107,17)(107,94.1077)
\drawline(107,17)(63,94.1077)
\drawline(107,17)(41,55.55385)
}

\allinethickness{1.6pt}
\drawline(128.6,30.2)(84.9,106.7)
\drawline(128.6,30.2)(41.4,80.8)
\drawline(128.6,30.2)(41.5,30)
\end{picture}
\caption{Initial seed in type $A_3$.}
\label{fig:type-a3-IN}
\end{figure}
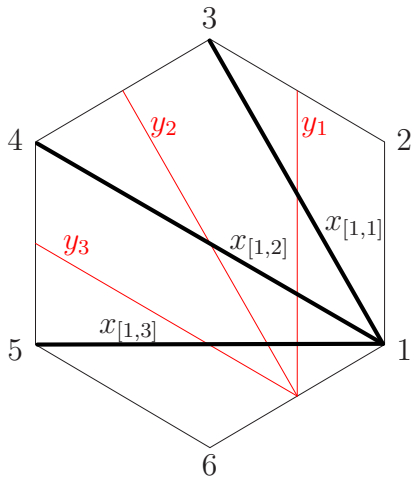

Clearly, the set of labels $\Pi(c)$ is in a bijection with the set of all
diagonals of $\Poly_{n+3}$ via the correspondence
\begin{equation}
\label{eq:two-parameterizations}
[i,j] \longleftrightarrow \langle i,j+2 \rangle
\quad (1 \leq i\leq j \leq n+1).
\end{equation}
Renaming each cluster variable $x_{[i,j]}$
in Theorem~\ref{th:ca-principal-minors-A-special}
by $x_{\langle i,j+2 \rangle}$, we see that, if we ignore the coefficients, then
the desired relations \eqref{eq:exchange-A-principal}
turn into \eqref{eq:ex-A-UC}.
Note that, under this relabeling, the initial cluster gets associated with the
triangulation of $\Poly_{n+3}$ by all the diagonals from the
vertex~$1$, see Figure~\ref{fig:type-a3-IN} (ignoring for the
moment the segments labeled by $y_1, y_2$ and $y_3$).

To complete the proof, it remains to verify the coefficients in the relations
\eqref{eq:exchange-A-principal}.
To do this, we use the geometric realization of the
universal coefficients (cf.~Proposition~\ref{pr:universal-ca-c} and \cite[Theorem~12.4]{ca4})
due to S.~Fomin and A.~Zelevinsky (the proof will appear elsewhere).

Consider the dual $(n+3)$-gon $\Poly'_{n+3}$ whose vertices are
the midpoints of the sides of $\Poly_{n+3}$.
For $i = 2, \ldots, n + 3$, we denote by $i'$ the midpoint of the side with vertices
$i-1$ and $i$; we also denote by $1'$ the midpoint of the side connecting $1$ and $n+3$.
We refer to the diagonals of $\Poly'_{n+3}$ as \emph{dual diagonals}.

\begin{proposition}[S.~Fomin, A.~Zelevinsky]
\label{pr:UC-type A}
After relabeling the cluster variables and the generators
of the coefficient semifield in the cluster algebra $\widetilde \Acal(c)$
in Proposition~\ref{pr:universal-ca-c} by the diagonals of $\Poly_{n+3}$ via
the correspondence \eqref{eq:two-parameterizations}, the coefficient $p^+_{ik,j\ell}$
(resp., $p^-_{ik,j\ell}$) in an exchange relation \eqref{eq:ex-A-UC}
turns into the product of the generators $p\langle a,b \rangle$
such that the dual diagonal $\langle a',b' \rangle$ is contained in the strip formed by $\langle i,j \rangle$
and $\langle k,\ell \rangle$ (resp., by $\langle i,\ell \rangle$
and $\langle j,k \rangle$).
\end{proposition}

\begin{example}
\label{ex:exchange-relation}
In type $A_3$, the exchange
relation between $x\langle 2,5 \rangle$ and
$x\langle 4,6 \rangle$ has the form
$$x\langle 2,5 \rangle \, x\langle 4,6 \rangle = p\langle 1,5 \rangle p\langle 2,5 \rangle \,
x\langle 2,4 \rangle \, x\langle 5,6 \rangle + p\langle 3,6 \rangle p\langle 4,6 \rangle \,
x\langle 2,6 \rangle \, x\langle 4,5 \rangle,$$
as illustrated by Figure~\ref{fig:type-a3-UC} (where we show only those
dual diagonals that contribute to the exchange relation);
note that by our convention,
$x\langle 4,5 \rangle = x\langle 5,6 \rangle = 1$.

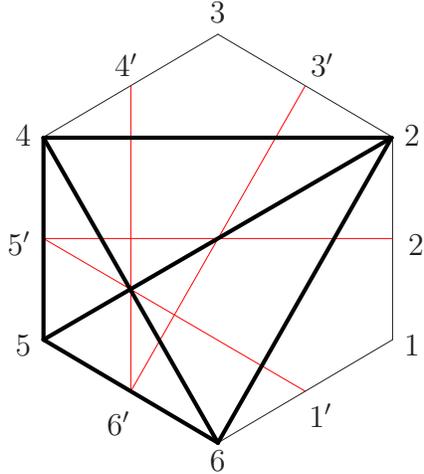
\begin{figure}[ht]
\setlength{\unitlength}{1.5pt}
\begin{picture}(160,120)(0,0)

\put(132,26){$1$}
\put(132,79){$2$}
\put(83,110){$3$}
\put(34,79){$4$}
\put(34,26){$5$}
\put(83,-3){$6$}

\put(108,8){$1'$}
\put(133,51.5){$2'$}
\put(108.4,96.4){$3'$}
\put(59,96.4){$4'$}
\put(32,51.5){$5'$}
\put(57,6){$6'$}

\drawline(85,4)(129,30)(129,81.1077)(85,107.1077)(41,81.1077)(41,30)(85,4)

\light{
\drawline(63,17)(63,94.1077)
\drawline(63,17)(107,94.1077)
\drawline(41,55.55385)(129,55.55385)
\drawline(41,55.55385)(107,17)
}

\allinethickness{1.6pt}
\drawline(41,30)(41,81)
\drawline(41,30)(85,4)
\drawline(41,30)(128.8,81)
\drawline(41,81)(85,4)
\drawline(41,81)(128.8,81)
\drawline(85,4)(128.8,81)
\end{picture}
\caption{Exchange relation between $x\langle 2,5 \rangle$ and
$x\langle 4,6 \rangle$.}
\label{fig:type-a3-UC}
\end{figure}

\end{example}

As shown in Section~\ref{sec:univcoeff}, the cluster algebra with
principal coefficients $\Acal(c)$ is obtained from $\widetilde \Acal(c)$
by the coefficient specialization
\eqref{eq:principal-specialization-of universal-coeffs}.
Using the relabeling in Proposition~\ref{pr:UC-type A}, this
specialization is described as follows: it sends
$p \langle 1, j+2 \rangle$ to $y_j$ for $j = 1, \dots, n$, and the rest of
the generators $p \langle i, j \rangle$ are sent to~$1$.
(This specialization is illustrated by Figure~\ref{fig:type-a3-IN}
showing the dual diagonals associated with $y_1, y_2, y_3$.)

Now consider any relation \eqref{eq:exchange-A-principal}.
Rewriting it in the form \eqref{eq:ex-A-UC}, we obtain
$$x_{\langle i,k+2 \rangle}\, x_{\langle j, \ell+2 \rangle}
 = p^+ x_{\langle i,j \rangle}\, x_{\langle k+2, \ell+2 \rangle}
+ p^- \, x_{\langle i,\ell+2 \rangle}\, x_{\langle j, k+2 \rangle} \ ,$$
where the coefficients $p^+$ and $p^-$ are obtained from those
given in Proposition~\ref{pr:UC-type A} by the above specialization.
By inspection of the corresponding dual diagonals, we conclude
that $p_+ = y_{j-1} y_j \cdots y_k$ and $p_- = 1$, verifying
the coefficients in \eqref{eq:exchange-A-principal} and finishing the proof of
Theorem~\ref{th:ca-principal-minors-A-special}.

\medskip

We conclude the paper by an example illustrating the
correspondence between various reduced double cells given by
\eqref{eq:c-tildec-isom-geometric}.

\begin{example}
Let $c = s_1 s_2$ and $\tilde c = s_2 s_1$ be the two Coxeter
elements in type~$A_2$, and let $L = L^{c,c^{-1}}$ and
$\widetilde L = L^{\tilde c,\tilde c^{-1}}$ be the corresponding
reduced double Bruhat cells in $G = SL_3(\CC)$.
The variety $L$ has been already described: it consists of matrices
$$M = \begin{pmatrix}v_1 & y_1 & 0 \\
                  1  &  v_2 & y_2 \\
                  0  & 1 & v_3
                  \end{pmatrix}$$
with nonzero $y_1, y_2, y_3$.
Using Propositions~\ref{pr:double-cell-by-eqs}
and~\ref{pr:Luv-equations}, the variety $\widetilde L$ can be
described as the set of matrices $\widetilde M = (m_{i,j}) \in G$
satisfying the conditions
\begin{itemize}
\item $\Delta_{[1,2],[2,3]}(\widetilde M) = \Delta_{[2,3],[1,2]}(\widetilde M)= 0$;
\item $m_{1,3} \neq 0$, $\Delta_{[1,2],\{1,3\}}(\widetilde M) \neq 0$;
\item $m_{3,1} = \Delta_{\{1,3\},[1,2]}(\widetilde M) = 1$.
\end{itemize}

The correspondence $L \to \widetilde L$ in \eqref{eq:c-tildec-isom-geometric} sends a
generic element
$$M = x_{-1}(u_1) x_{-2}(u_2)  x_{2}(t_2)
x_{1}(t_1) \in L$$
to $$\widetilde M = x_{1}(u_1) x_{-2}(u_2)  x_{2}(t_2)
x_{-1}(t_1^{-1}) \in \widetilde L.$$
Performing the matrix multiplication, we see that $M$ and
$\widetilde M$ are given by
$$M = \begin{pmatrix}u_1^{-1} & u_1^{-1} t_1 & 0 \\
                  1  &  u_1 u_2^{-1} + t_1 & u_1 u_2^{-1}t_2 \\
                  0  & 1 & u_2 + t_2
                  \end{pmatrix}, \quad
                  \widetilde M = \begin{pmatrix}
                  u_1 u_2^{-1} + t_1 & u_1 u_2^{-1}t_1^{-1}  & u_1 u_2^{-1}t_2 \\
                  u_2^{-1}  &  u_2^{-1} t_1^{-1} & u_2^{-1}t_2 \\
                  1  & t_1^{-1} & u_2 + t_2
                  \end{pmatrix}.$$
An easy calculation shows that the correspondence $M \mapsto
\widetilde M$ extends to an isomorphism $L \to \widetilde L$
given by
$$M = \begin{pmatrix}
                   v_1 & y_1 & 0 \\
                  1  &  v_2 & y_2 \\
                  0  & 1 & v_3
\end{pmatrix} \mapsto
\widetilde M = \begin{pmatrix}
                  v_2 & v_1 v_2y_1^{-1}-1 & y_2 \\
                  v_1 v_2-y_1  &  v_1 (v_1 v_2 y_1^{-1}-1) & v_1 y_2 \\
                  1  & v_1 y_1^{-1} & v_3
                  \end{pmatrix},$$
while the inverse isomorphism $\widetilde L \to  L$ is given by
$$\widetilde M = \begin{pmatrix}
                   m_{1,1} &  m_{1,1} m_{3,2}-1 &  m_{1,3} \\
                   m_{2,1}  &   m_{2,2} &  m_{2,3} \\
                  1  &  m_{3,2} &  m_{3,3}
\end{pmatrix} \mapsto$$
$$M = \begin{pmatrix}
                  m_{2,3}m_{1,3}^{-1} & m_{1,1} m_{2,3} m_{1,3}^{-1} - m_{2,1}  & 0 \\
                  1  &  m_{1,1}  &  m_{1,3}  \\
                  0  & 1 &  m_{3,3}
                  \end{pmatrix}.$$
\end{example}

\section*{Acknowledgments}

The authors benefited from the correspondence with Alexander Kirillov, Jr., and
discussions with Nathan Reading and David Speyer, concerning the
connections of their work to ours.
Part of this work was done during A.~Zelevinsky's stay at the
University of Bielefeld supported by a Humboldt Research Award; he
is grateful to Claus Ringel and other colleagues for their warm
hospitality and stimulating working conditions.


\begin{thebibliography}{xxx}

\bibitem{bfz}
A.~Berenstein, S.~Fomin, and A.~Zelevinsky, Parametrizations
of canonical bases and totally positive matrices, {\it Adv. in Math.}
{\bf 122} (1996), 49--149.

\bibitem{ca3}
A.~Berenstein, S.~Fomin, and A.~Zelevinsky,
Cluster algebras~III: Upper bounds and double Bruhat cells,
\textsl{Duke Math.~J.} \textbf{126} (2005), 1--52.

\bibitem{bz01}
A.~Berenstein and A.~Zelevinsky,
Tensor product multiplicities, canonical bases and totally positive varieties,
\textsl{Invent.\ Math.} \textbf{143} (2001), 77--128.

\bibitem{bourbaki} N.~Bourbaki,
{\em Lie Groups and Lie Algebras, Chapters 4-6},
Elements of Mathematics, Springer-Verlag, Berlin-Heidelberg-New
York, 2002.

\bibitem{fz-double}
S.~Fomin and A.~Zelevinsky, Double Bruhat cells and total
positivity, \textsl{J.~Amer.\ Math.\ Soc.} \textbf{12} (1999), 335--380.

\bibitem{fz-intel}
S.~Fomin and A.~Zelevinsky,
Total positivity: tests and parametrizations,
\textsl{Math. Intelligencer}   \textbf{22} (2000), no. 1, 23--33.

\bibitem{ca1}
S.~Fomin and A.~Zelevinsky,
Cluster algebras~I: Foundations,
\textsl{J.~Amer.\ Math.\ Soc.} \textbf{15} (2002), 497--529.

\bibitem{yga}
 S.~Fomin and A.~Zelevinsky,
$Y$-systems and generalized associahedra,
\textsl{Ann. in Math.} \textbf{158} (2003), 977--1018.

\bibitem{ca2}
S.~Fomin and A.~Zelevinsky,
Cluster algebras~II: Finite type classification,
\textsl{Invent.\ Math.} \textbf{154} (2003), 63--121.

\bibitem{ca4}
S.~Fomin and A.~Zelevinsky,
Cluster algebras~IV: Coefficients,
\textsl{Comp.\ Math.} \textbf{143} (2007), 112--164.

\bibitem{kirth}
A.~Kirillov Jr. and J.~Thind,
Coxeter elements and periodic Auslander-Reiten quiver,

\texttt{math/0703361}.

\bibitem{kostant-1} B.~Kostant,
The principal three-dimensional subgroups and the Betti numbers of a
complex simple Lie group, \textsl{Amer.\ J.\ Math.} \textbf{81}
(1959), 973--1032.

\bibitem{kostant-2} B.~Kostant,
The solution to a generalized Toda lattice and representation
theory, \textsl{Adv.\ in \ Math.} \textbf{34} (1979), no. 3, 195–-338.

\bibitem{reading-speyer} N.~Reading and D.~Speyer,
Cambrian fans, \texttt{math.CO/0606201}.

\bibitem{steinberg}
R.~Steinberg, Finite reflection groups,
\textsl{Trans.\ Amer.\ Math.\ Soc.} \textbf{91} (1959), 493--504.

\end{thebibliography}
\end{document}